\documentclass[final,10pt]{elsarticle}

\makeatletter
\def\ps@pprintTitle{%
	\let\@oddhead\@empty
	\let\@evenhead\@empty
	\let\@oddfoot\@empty
	\let\@evenfoot\@oddfoot
}
\makeatother

\usepackage{amssymb}
\usepackage{amsmath}
\usepackage{amsthm}

\usepackage{lineno}

\usepackage{amsfonts}
\usepackage{graphicx}
\usepackage{epstopdf}
\usepackage{bm,color,geometry,subfig}
\usepackage{eucal}
\usepackage{algpseudocode}
\usepackage{algorithm}
\usepackage{tabularx}
\usepackage{booktabs,mathrsfs}
\usepackage{multirow}

\newtheorem{theorem}{Theorem}[section]
\newtheorem{corollary}{Corollary}[theorem]
\theoremstyle{definition}
\newtheorem{definition}{Definition}[section]

\newcommand{\norv}{N}               
\newcommand{\svc}{\xi}        		
\newcommand{\pvc}{x}          		

\newcommand{\sv}{\bm{\svc}}    		
\newcommand{\pv}{\bm{\pvc}}    		
\newcommand{\pdf}{\rho}    	        
\newcommand{\pdfc}{\rho}            
\newcommand{\nstoch}{n_{\svc}}      
\newcommand{\nphy}{n_{\pvc}}        

\newcommand{\gpc}{\phi}             
\newcommand{\GPC}{\Phi}             
\newcommand{\pbf}{\psi}             

\newcommand{\ap}{{\scriptsize\mbox{ap}}}  %

\newcommand{\rd}{\mathrm{d}}        
\newcommand{\tol}{\mathtt{RelTol}}
\newcommand{\rTheta}{\mathrm{\Theta}}

\newcommand{\quasimatrix}[1]{\mathbb{#1}}
\newcommand{\trans}{\mathrm{T}}
\newcommand{\trunc}{p}
\newcommand{\coarse}{{\mathrm{c}}}
\newcommand{\samp}{{\mathrm{samp}}}
\newcommand{\fine}{{\mathrm{f}}}

\newcommand{\temp}{{\mathrm{temp}}}
\newcommand{\diag}{\mathrm{diag}}

\def\dsR{\mathbf{R}}

\begin{document}

\begin{frontmatter}



\title{An alternating low-rank projection approach for\\ partial differential equations with random inputs}

 \author[label1]{Guanjie Wang}\ead{wangguanjie0@126.com}
 \affiliation[label1]{organization={School of Statistics and Mathematics,
		Shanghai Lixin University of Accounting and Finance},
             city={Shanghai},
             postcode={201209},
             country={China}}
\author[label2]{Qifeng Liao\corref{cor1}}\ead{liaoqf@shanghaitech.edu.cn}
\affiliation[label2]{organization={School of Information Science and Technology, ShanghaiTech University},
             city={Shanghai},
             postcode={201210},
             country={China}}
\cortext[cor1]{Corresponding author}

\begin{abstract}
It is known that standard stochastic Galerkin methods face challenges when solving partial differential equations (PDEs) with random inputs. These challenges are typically attributed to the large number of required physical basis functions and stochastic basis functions. Therefore, it becomes crucial to select effective basis functions to properly reduce the dimensionality of both the physical and stochastic approximation spaces. In this study, our focus is on the stochastic Galerkin approximation associated with generalized polynomial chaos (gPC). We delve into the low-rank approximation of the quasimatrix, whose columns represent the coefficients in the gPC expansions of the solution. We conduct an investigation into the singular value decomposition (SVD) of this quasimatrix, proposing a strategy to identify the rank required for a desired accuracy. Subsequently, we introduce both a simultaneous low-rank projection approach and an alternating low-rank projection approach to compute the low-rank approximation of the solution for PDEs with random inputs. Numerical results demonstrate the efficiency of our proposed methods for both diffusion and Helmholtz problems.
\end{abstract}



\end{frontmatter}


\section{Introduction}

Over the past few decades, there has been rapid development in efficient numerical methods for solving partial differential equations (PDEs) with random inputs. This surge in interest has been driven by the necessity of performing uncertainty quantification when modeling practical problems, such as groundwater flow and acoustic scattering. The sources of uncertainty in these problems typically arise from lack of knowledge or accurate measurements of realistic model parameters, such as permeability coefficients and refraction coefficients~\cite{Xiu2002modeling,xiu2004two,Elman2005,David2009,Feng2015}.

Among uncertainty quantification approaches, this work specifically focuses on the stochastic Galerkin methods, which have been demonstrated to be effective in various disciplines~\cite{Ghanem2003,Xiu2010,bedelman16,David2016,jin2016well}. In the context of stochastic Galerkin methods, this work focus on the  generalized polynomial chaos (gPC) methods to discretize the stochastic parameter space~\cite{Xiu2002modeling,Xiu2002wiener}, and employ finite element methods to discretize the physical space~\cite{Elman2014}.

Designing efficient solvers is a crucial and challenging problem when dealing with the stochastic Galerkin methods, as it often leads to large coupled linear systems. Since the stochastic Galerkin methods typically result in linear systems formulated in the Kronecker form~\cite{Pellissetti2000,ullmann07,Powell2009}, various iterative algorithms that exploit the Kronecker structure of linear systems to reduce computational costs have been extensively studied in recent years~\cite{matthies2012solving,doostan2009least,powell2015,LeeElman16,elman2018low,lee2019low,wang2024reduced}. Additionally, it has been shown that the desired stochastic Galerkin solution can be approximated through low-rank approximation, which further reduces computational costs~\cite{benner2015low}.

In this work, we delve into low-rank approximations for solving PDEs with random inputs, and propose an alternating low-rank projection (AltLRP) approach.  To achieve this, we first introduce the concept and properties of the~\textit{quasimatrix}~\cite{battles2005numerical,Townsend2014computing,townsend2015continuous}. Then, we investigate the singular value decomposition (SVD) of quasimatrices associated with the truncated gPC expansions of the solutions. The relationship between the singular values of these  quasimatrices is explored in detail, accompanied by a thorough theoretical analysis. Building upon these findings, a systematic strategy is proposed to identify the rank for a desired accuracy, and  a simultaneous low-rank projection (SimLRP) approach is developed, where reduced basis functions of the stochastic and physical spaces are constructed simultaneously. Based on SimLRP, our AltLRP approach is proposed to refine the solution, with the outcome of SimLRP serving as the initial guess for AltLRP. The effectiveness of AltLRP is illustrated through numerical studies.  

We note that there are various methods developed to exploit the inherent low-rank structure of solutions to stochastic PDEs. These include approaches that incorporate tensor-based representations of the input data and operators~\cite{espig2014efficient,dolgov2015polynomial}, as well as low-rank iterative solvers for the linear systems arising from stochastic Galerkin methods~\cite{matthies2012solving,doostan2009least,elman2018low,lee2019low}. In the former case, a tensor train format is constructed for the coefficient random field involved in the stochastic PDEs, enabling efficient assembly of the stochastic Galerkin matrix in the tensor train format with the same ranks~\cite{dolgov2015polynomial}. In this work, we adopt a truncated Karhunen–Lo{\'e}ve (KL) expansion~\cite{Ghanem2003,Elman2007} to represent the random input field and focus on the standard stochastic Galerkin method, without relying on tensorized representations of the coefficient or operator.

A key distinction between the proposed approach and existing low-rank iterative solvers lies in how the rank is handled throughout the computation. Most existing solvers operate directly on the full linear systems arising from stochastic Galerkin methods and apply low-rank techniques during the iterations. These methods typically either introduce truncation operators at each iteration to control rank growth~\cite{matthies2012solving,LeeElman16}, or adaptively increase the rank to achieve a desired level of accuracy~\cite{doostan2009least,lee2019low}. In contrast, the methods proposed in this work determine the rank by analyzing the singular values of quasimatrices associated with the truncated gPC expansion of the solution on a coarse physical grid. This pre-analysis enables the construction of a fixed-rank approximation, thereby avoiding the need for dynamic truncation or rank adaptation during iterative procedures for solving the linear systems, such as preconditioned conjugate gradient (PCG).

Another significant difference lies in the formulation of the linear system. Instead of solving the original linear system directly, we first construct a much smaller low-rank projected linear system based on the identified low-rank subspace. Once this reduced system is solved, an approximation to the original solution can be efficiently reconstructed. As a result, the proposed method significantly reduces the computational cost compared to existing approaches.

The outline of this work is as follows. In Section~\ref{sec:probset}, we describe the PDE formulation considered in this study and introduce the stochastic Galerkin method. In Section~\ref{sec:altlrp}, our main theoretical analysis concerning the relationship between the singular values of quasimatrices 
associated with gPC is discussed, and our AltLRP approach is presented. Numerical results are discussed in Section~\ref{sec:numericaltests}, and Section~\ref{sec:concusion} concludes the paper.

\section{Problem setting and stochastic Galerkin method}\label{sec:probset}
This section  describes the problems considered in this study, revisits the stochastic Galerkin method, and  presents the resulting linear system from the stochastic Galerkin method.

\subsection{Problem setting}
Let $D$ be a bounded and connected physical domain in $\mathbf{R}^d~(d = 2, 3)$ with a polygonal boundary denoted as $\partial D$ and $\pv \in \mathbf{R}^d$ represent a physical variable. Let $\sv = (\svc_1,\ldots,\svc_{\norv})$ be a random vector of dimension $\norv$, where the image of $\svc_i$ is denoted by $\Gamma_i$, and the probability density function of~$\svc_i$  is denoted by $\rho_{i}(\svc_i)$. We further assume that the components of $\sv$, i.e., $\svc_1,\ldots,\svc_{\norv}$ are mutually independent, the image of $\sv$ is then given by $\Gamma = \Gamma_1\times\cdots\times\Gamma_{\norv}$, and the probability density function of $\sv$ is given by $\pdf(\sv) = \prod_{i=1}^{\norv}\pdfc_{i}(\svc_i)$. In this work, we consider the following partial differential equations (PDEs) with random inputs, i.e.,
\begin{equation}\label{eq:spde}
\begin{cases}
\mathcal{L}(\pv,\sv,u(\pv,\sv)) = f(\pv)\ & \forall~ \pv\in D,\\
\mathcal{B}(\pv,\sv,u(\pv,\sv)) = g(\pv)\ &\forall~ \pv\in\partial D,
\end{cases}
\end{equation}
where $\mathcal{L}$ is a linear partial differential operator with respect to the physical variable, and $\mathcal{B}$ is a boundary operator. Both  operators may have random coefficients. The source function is denoted by $f(\pv)$, and $g(\pv)$ specifies the boundary conditions.  Additionally, we assume that $\mathcal{L}$ and $\mathcal{B}$ are affinely dependent on the random inputs,~i.e., 
\begin{eqnarray}
\mathcal{L}(\pv,\sv,u(\pv,\sv)) = \sum_{i=1}^{K}\rTheta_{\mathcal{L}}^{(i)}(\sv)\mathcal{L}_i(\pv,u(\pv,\sv)), \label{eq:affine1}\\ 
\mathcal{B}(\pv,\sv,u(\pv,\sv)) = \sum_{i=1}^{K}\rTheta_{\mathcal{B}}^{(i)}(\sv)\mathcal{B}_i(\pv,u(\pv,\sv)),\label{eq:affine2}
\end{eqnarray} 
where $\{\mathcal{L}_{i}\}_{i=1}^{K}$ are parameter-independent  linear differential operators, and $\{\mathcal{B}_{i}\}_{i=1}^{K}$ are parameter-independent  boundary operators. Both $\rTheta_{\mathcal{L}}^{(i)}(\sv)$ and $\rTheta_{\mathcal{B}}^{(i)}(\sv)$ take values in $\dsR$ for $i=1,\ldots,K$.

It is of interest to design a surrogate model for the problem~\eqref{eq:spde} or calculate statistics of the stochastic solution $u(\pv,\sv)$, such as the mean and the variance. For simplicity of the presentation, we consider problems that satisfy homogeneous Dirichlet boundary conditions. However, it is noteworthy that the approach we present can be readily extended to other arbitrary well-posed boundary conditions.

\subsection{Stochastic Galerkin method}\label{sec:SGM}

We provide a concise overview of the stochastic Galerkin method, as detailed in~\cite{Ghanem2003,Babuska2004}. To facilitate our discussion, we begin by introducing the variational formulation of~\eqref{eq:spde}. In order to do this, it is necessary to establish some key notations. Firstly, we define the Hilbert spaces $L^2(D)$ and $L^2_{\pdf}(\Gamma)$ as follows:
\begin{eqnarray*}
	L^2(D)&:=&\left\{\pbf(\pv): D \to \dsR\ \bigg|\ \int_D \pbf^*(\pv)\pbf(\pv)\, \rd\pv <\infty \right\}, \\
	L_{\pdf}^2(\Gamma)&:=&\left\{\gpc(\sv): \Gamma \to \dsR\ \bigg|\ \int_{\Gamma}\pdf(\sv) \gpc^*(\sv)\gpc(\sv)\, \rd\sv <\infty \right\},
\end{eqnarray*} 
which are equipped with the inner products
\begin{eqnarray*}\label{eq:inner1}
	\langle \pbf(\pv),\widetilde{\pbf}(\pv)\rangle_{L^2}:=\int_D \widetilde{\pbf}^*(\pv)\pbf(\pv)\,\rd\pv, \
	\langle \gpc(\sv),\widetilde{\gpc}(\sv)\rangle_{L^2_\pdf}:=\int_{\Gamma}\pdf(\sv) \widetilde{\gpc}^*(\sv)\gpc(\sv)\,\rd\sv.\label{eq:inner2}
\end{eqnarray*}
Following~\cite{Babuska2004}, we define the tensor space of $L^2(D)$ and $L^2_{\pdf}(\Gamma)$ as
\begin{equation}\notag
L^2(D)\otimes L^2_{\pdf}(\Gamma):=\left\{ \omega(\pv,\sv) \bigg| \omega(\pv,\sv)=\sum_{i=1}^{n}\pbf_i(\pv)\gpc_i(\sv), \pbf_i(\pv)\in L^2(D),\, \gpc_i(\sv)\in L^2_{\pdf}(\Gamma), n\in \mathbf{N}^+ \right\},
\end{equation}
which is equipped with the inner product
\begin{equation}\notag
\left\langle \omega(\pv,\sv),\widetilde{\omega}(\pv,\sv)\right\rangle_{L^2\otimes L_{\pdf}^2} =\int_\Gamma\int_{D}\pdf(\sv)\widetilde{\omega}^*(\pv,\sv) \omega(\pv,\sv) \rd\pv \rd\sv.
\end{equation}
In general, the variational form of~\eqref{eq:spde} can be expressed as follows: find $u$ in $W = X(D)\otimes L^2_{\pdf}(\Gamma )$ such that \begin{equation}\label{eq:weak}
	\mathcal{A}(u,\omega) = \mathcal{F}(\omega),\ \forall~ \omega\in W,
	\end{equation}
	where $X(D)$ is an appropriate subspace of $L^2(D)$ associated with $D$, and 
	\begin{equation}\notag
	\mathcal{A}(u,\omega):=\left\langle \mathcal{L}(\pv,\sv,u(\pv,\sv)),\omega(\pv,\sv)\right\rangle_{L^2\otimes L^2_{\pdf}}, 
	\ \mathcal{F}(\omega):= \left\langle f(\pv),\omega(\pv, \sv) \right\rangle_{L^2 \otimes L^2_{\pdf}}.
	\end{equation}
The stochastic Galerkin method aims to approximate~\eqref{eq:weak} within a finite-dimensional subspace of~$W$. This subspace is defined as
\begin{equation}\notag
W_h^p:= X_h \otimes S_p :=\mathrm{span} \left\{\pbf_s(\pv)\gpc_j(\sv)\left| \pbf_s\in X_h, \gpc_j\in S_p \right.\right\},
\end{equation}
where   
\begin{equation}\label{eq:phy_stoch_full}
X_h = \mathrm{span}\left\{\pbf_s(\pv)\right\}_{s=1}^{\nphy}\subseteq X(D), \ 
S_p= \mathrm{span}\left\{\gpc_j(\sv)\right\}_{j=1}^{\nstoch}\subseteq L^2_{\pdf}(\Gamma).
\end{equation}
In this context, $\{\pbf_s(\pv)\}_{s=1}^{\nphy}$ represents physical basis functions, while $\{\gpc_j(\sv)\}_{j=1}^{\nstoch}$ represents stochastic basis functions. 

In this study, we employ quadrilateral elements and piecewise bilinear basis functions to discretize the physical domain, while for the stochastic domain, we utilize the generalized polynomial chaos basis functions~\cite{Xiu2002modeling}. Note that in this situation,~$\{\gpc_j(\sv)\}_{j=1}^{\nstoch}$ is orthonormal in $L^2_{\rho}(\Gamma)$.

Assuming that $u^{\ap}(\pv,\sv)$ is defined as
\begin{equation}\label{eq:appsg}
u^{\ap}(\pv,\sv):= \sum_{j=1}^{\nstoch}\sum_{s=1}^{\nphy}u_{sj}\pbf_s(\pv)\gpc_j(\sv).
\end{equation}
Since $\mathcal{L}$ is affinely dependent on the random inputs (see~\eqref{eq:affine1}), we substitute~\eqref{eq:appsg} into~\eqref{eq:weak} to obtain~\cite{Powell2009}, 
\begin{equation}\label{eq:linsg}
\left(\sum_{i=1}^{K}\bm{G}_i\otimes \bm{A}_i\right)\bm{u} = \bm{h}\otimes\bm{f}, 
\end{equation}
where $\{\bm{G}_i\}_{i=1}^K$ are matrices of size $\nstoch\times \nstoch$, and $\bm{h}$ is a column vector of length $\nstoch$. They are defined via  
\begin{equation}\label{eq:stoch_mat}
\bm{G}_i(j,k) = \langle \rTheta_{\mathcal{L}}^{(i)}(\sv)\gpc_j(\sv),\gpc_k(\sv) \rangle_{L^2_\pdf}, \ \bm{h}(k) = \langle 1,\gpc_k(\sv) \rangle_{L^2_\pdf}.
\end{equation}
The matrices $\bm{A}_i$ and the vector $\bm{f}$ in~\eqref{eq:linsg} are defined through
\begin{equation}\label{eq:phy_mat}
\bm{A}_i(s,t) = \left\langle\mathcal{L}_i\pbf_s,\pbf_t \right\rangle_{L^2},\ \bm{f}(s) = \left\langle f,\pbf_s \right\rangle_{L^2}, 
\end{equation}
where $s = 1,\ldots\nphy$ and $t = 1,\ldots, \nphy$. The vector $\bm{u}$ in~\eqref{eq:linsg} is a column vector of length~$\nphy\times\nstoch$, and is defined by
\begin{equation}\notag
\bm{u}=\begin{bmatrix}
\bm{u}_1 \\
\vdots\\
\bm{u}_{\nstoch}
\end{bmatrix}, \ \mbox{where}\ 
\bm{u}_j=\begin{bmatrix}
{u}_{1j} \\
\vdots\\
{u}_{\nphy j}
\end{bmatrix}, \ j= 1,\ldots, \nstoch.
\end{equation}
The linear system~\eqref{eq:linsg} can be rewritten in a matrix form
\begin{equation}\label{eq:linsg_matrix}
\sum_{i=1}^{K} \bm{A}_i \bm{U}\bm{G}_i^\trans = \bm{f} \bm{h}^\trans,
\end{equation}
where 
$\bm{U} = [\bm{u}_1,\ldots,\bm{u}_{\nstoch}]$.

\section{Alternating low-rank projection approach}\label{sec:altlrp}
In this section, we introduce an alternating low-rank projection (AltLRP) approach for PDEs with random inputs. To this end, we first present the low-rank approximation of a quasimatrix and investigate the singular value decomposition~(SVD) of the quasimatrix. The columns of this quasimatrix represent the coefficients of the truncated generalized polynomial chaos expansions corresponding to the solution for~\eqref{eq:spde}. Additionally, we explore the SVD of another quasimatrix, where the columns represent the values of the truncated generalized polynomial chaos expansions at random samples of random inputs~(i.e.,~$\sv$). The relationship between the singular values of these two quasimatrices is also examined.

Based on these findings, we propose a simultaneous low-rank projection (SimLRP) approach for PDEs with random inputs, in which the reduced basis functions of the stochastic and physical spaces are constructed simultaneously. To further improve the performance, we introduce an alternating low-rank projection (AltLRP) approach,  in which low-rank projections in the stochastic and physical spaces are applied alternately. AltLRP capitalizes on the outcomes of SimLRP by employing them as an initial guess, leveraging pre-established insights to facilitate a refined and robust convergence toward accurate solutions.

\subsection{The low-rank approximation of a quasimatrix}
Let us consider the generalized polynomial chaos expansion of $u(\pv,\sv)$, i.e., the solution to~\eqref{eq:spde}. Suppose that the generalized polynomial chaos expansion of $u(\pv,\sv)$ is given by
\begin{equation}\notag
u(\pv,\sv) = \sum_{j=1}^{\infty} u_j(\pv)\gpc_{j}(\sv),
\end{equation}
and we truncate this expansion at a gPC order $p$, i.e.,
\begin{equation}\label{eq:trunc_gpc}
u^{\trunc}(\pv,\sv) = \sum_{j=1}^{\nstoch} u_j(\pv)\gpc_{j}(\sv),\ \ (\pv,\sv) \in D\times\Gamma,
\end{equation}
where $\nstoch = (N+p)!/(N!p!)$. Rewrite~\eqref{eq:trunc_gpc} in an analogues matrix form
\begin{equation}\label{eq:trunc_matrix}
u^{\trunc}(\pv,\sv) = \quasimatrix{U}^{\trunc}_{D\times\nstoch}(\quasimatrix{S}^{\trunc}_{\Gamma\times \nstoch})^*, 
\end{equation}
where
\begin{equation}\notag
\quasimatrix{U}^{\trunc}_{D\times\nstoch} = \left[\begin{array}{llll}
\Big|    	&   \Big|  	&        	    &  \Big|        \\ 
u_1 		& u_2 		& \ldots 		& u_{\nstoch} \\
\Big|    	&   \Big|  	&        	    &  \Big|      
\end{array}\right],  \ \quasimatrix{S}^{\trunc}_{\Gamma \times\nstoch} = \left[\begin{array}{llll}
\Big|    	&   \Big|  	&        	    &  \Big|        \\ 
\gpc_1^* 		& \gpc_2^* 		& \ldots 		& \gpc_{\nstoch}^* \\
\Big|    	&   \Big|  	&        	    &  \Big|      
\end{array}\right].
\end{equation}
We refer to $\quasimatrix{U}^{\trunc}_{D\times\nstoch}$ and 
$\quasimatrix{S}^{\trunc}_{\Gamma \times\nstoch}$ as $D\times \nstoch$ and $\Gamma \times\nstoch$ quasimatrices, respectively, since one index of the rectangular matrix becomes continuous while the other remains discrete~\cite{townsend2015continuous}. We use the superscript $\trans$ to denote the transpose and $*$ to denote the conjugate transpose of a quasimatrix. More precisely, suppose that
\begin{equation}\notag
\quasimatrix{A} = \left[\begin{array}{llll}
\Big|    	&   \Big|  	&        	    &  \Big|        \\ 
a_1 		& a_2 		& \ldots 		& a_n \\
\Big|    	&   \Big|  	&        	    &  \Big| 
\end{array}\right]
\end{equation}
is a $D\times n$ quasimatrix, then 
\begin{equation}\notag
\quasimatrix{A}^\trans = \left[\begin{array}{ccc}
\text{---} 	    &  a_1	&  \text{---}           \\ 
\text{---}      &  a_2 		&  \text{---}\\
&   \vdots  	&       \\      
\text{---}  	&  a_n & \text{---}
\end{array}\right],\ 
\quasimatrix{A}^* = \left[\begin{array}{ccc}
\text{---} 	    &  a_1^*	&  \text{---}           \\ 
\text{---}      &  a_2^* 		&  \text{---}\\
&   \vdots  	&       \\      
\text{---}  	&  a_n^* & \text{---}
\end{array}\right].
\end{equation}
Thus $(\quasimatrix{S}^{\trunc}_{\Gamma\times \nstoch})^*$ in~\eqref{eq:trunc_matrix} is given by
\begin{equation}\notag
(\quasimatrix{S}^{\trunc}_{\Gamma\times \nstoch})^*=
\left[\begin{array}{ccc}
\text{---} 	    &  \gpc_1	&  \text{---}           \\ 
\text{---}      &  \gpc_2 		&  \text{---}\\
&   \vdots  	&       \\      
\text{---}  	&  \gpc_{\nstoch} & \text{---}
\end{array}\right].
\end{equation}
Note that  if $a_i (i=1,\ldots,n)$ are real functions, then $\quasimatrix{A}^\trans = \quasimatrix{A}^*$. The product of~$\quasimatrix{A}^*$ and $\quasimatrix{A}$ is defined as
\begin{equation}\notag
\quasimatrix{A}^*\quasimatrix{A} = \left[\begin{array}{cccc}
a_1^*a_1 	    &  a_1^*a_2	&  \ldots & a_1^*a_n \\ 
a_2^*a_1        &  a_2^*a_2 &  \ldots &  a_2^*a_n\\
\vdots	        &   \vdots  &  \ldots & \vdots   \\      
a_n^*a_1        &  a_n^*a_2 &  \ldots &  a_n^*a_n\\
\end{array}\right].
\end{equation}
For a matrix, QR factorization exists. The QR factorization of a $D\times n$ quasimatrix~$\quasimatrix{A}$ is a straightforward extension of the matrix case. Following~\cite{townsend2015continuous}, we first give the definition of the QR factorization of a quasimatrix.
\begin{definition}[QR factorization of a quasimatrix]
	Let $\quasimatrix{A}$ be a $D\times n$ quasimatrix. A QR factorization of $\quasimatrix{A}$ is a factorization $\quasimatrix{A} = \quasimatrix{Q}\bm{R}$,
	where $\quasimatrix{Q}$ is a $D \times n$ quasimatrix with orthonormal columns and $\bm{R}$ is an $n \times n$ upper triangular
	matrix with non-negative real numbers on the diagonal.
\end{definition}
The QR factorization of a $D\times n$ quasimatrix has the following fundamental properties~\cite{battles2005numerical,Townsend2014computing,townsend2015continuous}.
\begin{theorem}
	Every $D\times n$ quasimatrix has a QR factorization, which can be calculated by Gram–Schmidt orthogonalization. If the columns of $\quasimatrix{A}$ are linearly independent, the QR factorization is unique.
\end{theorem}

\begin{proof}
	The existence is evident through the Gram–Schmidt process, and we focus here on providing the proof of uniqueness\footnote{The proof is available in~\cite{battles2005numerical}; however, we are unable to access the literature in electronic or hard copy form.}. Suppose that
	\begin{equation}\notag
	\quasimatrix{A}=\quasimatrix{Q}_1\bm{R}_1 = \quasimatrix{Q}_2\bm{R}_2.
	\end{equation}
	Since the columns of $\quasimatrix{A}$ are linearly independent, $\bm{R}_1$ and $\bm{R}_2$ are invertible, and the diagonal elements of $\bm{R}_1$ and $\bm{R}_2$ are positive.
	Thus we have 
	\begin{equation}\notag
	\quasimatrix{Q}_1 = \quasimatrix{Q}_2\bm{R}_2\bm{R}_1^{-1},
	\end{equation}
	which implies
	\begin{equation}\notag
	\quasimatrix{Q}_1^*\quasimatrix{Q}_1 = (\bm{R}_2\bm{R}_1^{-1})^{*}\quasimatrix{Q}_2^{*}\quasimatrix{Q}_2\bm{R}_2\bm{R}_1^{-1}.
	\end{equation}
	Since the columns of $\quasimatrix{Q}_1$ and $\quasimatrix{Q}_2$ are orthonormal, by integration over $D$, we obtain 
	\begin{equation}\label{eq:qr_proof1}
	\bm{E} = (\bm{R}_2\bm{R}_1^{-1})^*\bm{R}_2\bm{R}_1^{-1}
	\end{equation}
	Note that $\bm{R}_1$ and $\bm{R}_2$ are upper triangular
	matrix, and thus $\bm{R}^{-1}$ and $\bm{R}_2\bm{R}_1^{-1}$ are also upper triangular
	matrices. The diagonal elements of $\bm{R}_2\bm{R}_1^{-1}$ are given by $\bm{R}_2(i,i)/\bm{R}_1(i,i)$.
	
	Suppose that
	\begin{equation}\notag
	\bm{R}_2\bm{R}_1^{-1} = 
	\begin{bmatrix}
	r_{11} & r_{12} & \ldots & r_{1n}\\
	0 & r_{22} & \ldots & r_{2n}\\
	0 &  0 &  \ddots       & \vdots\\
	0 &  0 &  \ldots       & r_{nn} 
	\end{bmatrix},\  (\bm{R}_2\bm{R}_1^{-1})^{*} = 
	\begin{bmatrix}
	r_{11}^* &  0 & \ldots &  0\\
	r_{12}^* & r_{22}^* & \ldots & 0\\
	\vdots &  \vdots &  \ddots       & \vdots\\
	r_{1n}^* &  r_{2n}^* &  \ldots       & r_{nn}^* 
	\end{bmatrix}.
	\end{equation}
	Since $r_{ii}=\bm{R}_2(i,i)/\bm{R}_1(i,i)>0$, by~\eqref{eq:qr_proof1}, we have
	$\bm{R}_2\bm{R}_1^{-1}=\bm{E}$, and thus 
	\begin{displaymath}
	\bm{R}_1=\bm{R}_2,\ \quasimatrix{Q}_1 = \quasimatrix{Q}_2.
	\end{displaymath}
\end{proof}

\begin{definition}[SVD decomposition of a quasimatrix]
	Let $\quasimatrix{A}$ be a $D\times n$ quasimatrix. An SVD of $\quasimatrix{A}$ is a factorization $ \quasimatrix{A} = \quasimatrix{U} \bm{\Sigma V}^*$, where $\quasimatrix{U}$
	is a $D\times n$ quasimatrix with orthonormal columns, $\bm{\Sigma}$ is an $n\times n$ diagonal matrix with diagonal
	entries $\sigma_1\geq \sigma_2\geq \ldots \sigma_n\geq 0 $ and $\bm{V}$ is an~$n\times n$ unitary matrix.
\end{definition}

\begin{theorem}\label{thm:quasi_svd}
	Every $D\times n$ quasimatrix has an SVD, which can be calculated by computing a QR decomposition $\quasimatrix{A} = \quasimatrix{Q}\bm{R}$ followed by a matrix SVD of the triangular factor, $\bm{R} = \bm{U}_1\bm{\Sigma V}^*$; an SVD of $\quasimatrix{A}$ is
	then obtained as $(\quasimatrix{Q}\bm{U}_1)\bm{\Sigma V}^*$. The singular values are unique, and the singular vectors corresponding to
	simple singular values are also unique up to complex signs. The
	rank $r$ of $\quasimatrix{A}$ is the number of nonzero singular values.
\end{theorem}

\begin{proof}
	The existence is obvious, and the proof of other conclusions can be found in~\cite{Townsend2014computing}. 
\end{proof}

\begin{theorem}\label{th:snap_equv}
	Suppose that $\sv_1,\ldots\sv_n$ are random samplings from $\sv$, and the quasimatrix $\quasimatrix{U}^{\samp}_{D\times n}$ is defined as
	\begin{equation}\notag
	\quasimatrix{U}^{\samp}_{D\times n} = \left[\begin{array}{llll}
	\Big|    	&   \Big|  	&        	    &  \Big|        \\ 
	u^{\trunc}(\pv,\sv_1) 		& u^{\trunc}(\pv,\sv_2) 		& \ldots 		& u^{\trunc}(\pv,\sv_n) \\
	\Big|    	&   \Big|  	&        	    &  \Big|      
	\end{array}\right].
	\end{equation}
	Then the number of non-zero singular values of $\quasimatrix{U}^{\samp}_{D\times n}$ is less than or equal to $\nstoch$. Denote the first $\nstoch$ singular values of $\quasimatrix{U}^{\samp}_{D\times n}$ by $\lambda_1\geq\lambda_2\geq\ldots\geq\lambda_{\nstoch}\geq0$, then for large $n$ we have
	\begin{equation}\notag
	\sqrt{\sum_{j=1}^{\nstoch}\left|\sigma_j^2-\frac{\lambda_j^2}{n}\right|^2}\left/\sum_{j=1}^{\nstoch}\sigma_j^2\right.\lesssim c\frac{\nstoch}{\sqrt{n}}\gamma.
	\end{equation}
	where $\sigma_1\geq\sigma_2\geq\ldots\geq\sigma_{\nstoch}\geq0$ are the singular values of $\quasimatrix{U}^{\trunc}_{D\times \nstoch}$, $c$ is a constant and~$\gamma$ is a random variable.
	
\end{theorem}

\begin{proof}
	Suppose that the SVD decomposition of $\quasimatrix{U}^{\trunc}_{D\times \nstoch}$ is given by 
	\begin{equation}\notag
	\quasimatrix{U}^{\trunc}_{D\times \nstoch} = \quasimatrix{U}\bm{\Sigma}\bm{V}^*,
	\end{equation}
	then we have 
	\begin{equation}\notag
	\quasimatrix{U}^{\samp}_{D\times n} = \quasimatrix{U}\bm{\Sigma}\bm{V}^*\bm{\GPC}^*,
	\end{equation}
	where
		\begin{equation}\notag
		\quasimatrix{U}^{\samp}_{D\times n} = \left[\begin{array}{llll}
		\Big|    	&   \Big|  	&        	    &  \Big|        \\ 
		u^{\trunc}(\pv,\sv_1) 		& u^{\trunc}(\pv,\sv_2) 		& \ldots 		& u^{\trunc}(\pv,\sv_n) \\
		\Big|    	&   \Big|  	&        	    &  \Big|      
		\end{array}\right], 
		\bm{\GPC}^* = \left[\begin{array}{cccc}
		\gpc_1(\sv_1) 	    		&  \gpc_1(\sv_2)			&  \ldots & \gpc_{1}(\sv_n) \\ 
		\gpc_2(\sv_1)      			&  \gpc_2(\sv_2) 			&  \ldots &  \gpc_{2}(\sv_n)\\
		\vdots	            		&   \vdots  				&  \ldots & \vdots   \\      
		\gpc_{\nstoch}(\sv_1)      	&  \gpc_{\nstoch}(\sv_n) 	&  \ldots &  \gpc_{\nstoch}(\sv_n)\\
		\end{array}\right].
		\end{equation}
	Since the rank of $\quasimatrix{U}^{\samp}_{D\times n}$ is less than or equal to $\nstoch$, the number of  non-zero singular values is less than or equal to $\nstoch$. The non-zero singular values of $\quasimatrix{U}^{\samp}_{D\times n}$ can be given by the square root of the eigenvalues of ~\cite{trefethen2022numerical}
	\begin{equation}
	\bm{T} := 
	\bm{\Sigma}\bm{V}^*\bm{\GPC}^*\bm{\GPC}\bm{V}\bm{\Sigma}^*.
	\end{equation}
	Now let us consider the eigenvalues of $\bm{T}$. 
	Note that $\{\gpc_1,\ldots,\gpc_{\nstoch}\}$ is orthonormal, that is
	\begin{equation}\notag
	\left\langle \gpc_{i}, \gpc_{j} \right\rangle_{L^2_{\pdf}} = \int_{\Gamma}\gpc_{j}^*(\sv)\gpc_{i}(\sv)\rd\sv = \delta_{ij},
	\end{equation}
	where $\delta_{ij}$ is the Kronecker's delta function. Thus for large $n$, we have~\cite{Caflisch1998}
	\begin{equation}\notag
	\frac{1}{n}\sum_{k=1}^{n}\gpc_{i}(\sv_k)\gpc^*_{j}(\sv_k)-\delta_{ij} \approx c_{ij}\frac{1}{\sqrt{n}}\gamma_{ij},    \ i,j = 1\ldots,\nstoch,
	\end{equation}
	where $c_{ij}$ is a constant, and $\gamma_{ij}$ is a standard normal random variable. This implies that	
	\begin{equation}\notag
	\left\Vert\frac{1}{n}\bm{\GPC}^*\bm{\GPC}-\bm{E}\right\Vert_{F}\lesssim c\frac{\nstoch}{\sqrt{n}} \gamma,
	\end{equation}
	where $c=\max\{c_{ij}\}$,  is  a constant, and $\gamma=\max\{|\gamma_{ij}|\}$ is a  random variable. Suppose that 
	\begin{equation}\notag
	\frac{1}{n}\bm{\GPC}^*\bm{\GPC} = \bm{E} +\bm{\Delta}, \ \mbox{where}\ \Vert\bm{\Delta}\Vert_{F}\lesssim  c\frac{\nstoch}{\sqrt{n}} \gamma.
	\end{equation}
	Note that
	\begin{equation}\notag
	\begin{split}
	\bm{T}
	=&(\sqrt{n}\bm{\Sigma})\bm{V}^*\left(\frac{1}{n}\bm{\GPC}^*\bm{\GPC}\right)\bm{V}(\sqrt{n}\bm{\Sigma})
	= (\sqrt{n}\bm{\Sigma}) \bm{V}^*(\bm{E}+\bm{\Delta})\bm{V}(\sqrt{n}\bm{\Sigma})\\
	=& n\bm{\Sigma}^2 + n(\bm{V}\bm{\Sigma})^*\bm{\Delta}(\bm{V}\bm{\Sigma}).
	\end{split}
	\end{equation}
	It is clear that the eigenvalues of $n\bm{\Sigma}^2$  are $n\sigma_1^2,\ldots,n\sigma_{\nstoch}^2$, and the eigenvalues of $\bm{T}$ are $\lambda_1^2,\ldots,\lambda_{\nstoch}^2$. Since $\bm{T}$ is a normal matrix, we have~\cite{hoffman2003variation,bhatia2007perturbation,li2006eigenvalue}
	\begin{equation}\notag
	\sqrt{\sum_{j=1}^{\nstoch}|n\sigma_j^2-\lambda_j^2|^2}\leq \Vert  n(\bm{V}\bm{\Sigma})^*\bm{\Delta}(\bm{V}\bm{\Sigma}) \Vert_{F}.
	\end{equation}
	Note that
	\begin{equation}\notag
	\begin{split}
	\Vert  n(\bm{V}\bm{\Sigma})^*\bm{\Delta}(\bm{V}\bm{\Sigma}) \Vert_{F}\leq n\Vert\Delta\Vert_{F}\Vert \bm{V}\bm{\Sigma}\Vert_{F}^2=n\Vert\bm{\Sigma}\Vert_{F}^2\Vert\bm{\Delta}\Vert_{F}
	\lesssim c\sqrt{n}\nstoch\Vert\bm{\Sigma}\Vert_{F}^2\gamma,
	\end{split}
	\end{equation}
	which implies 
	\begin{equation}\notag
	\sqrt{\sum_{j=1}^{\nstoch}|n\sigma_j^2-\lambda_j^2|^2}\lesssim c\sqrt{n}\nstoch\Vert\bm{\Sigma}\Vert_{F}^2\gamma.
	\end{equation}
	This further implies 
	\begin{equation}\notag
	\sqrt{\sum_{j=1}^{\nstoch}\left|\sigma_j^2-\frac{\lambda_j^2}{n}\right|^2}\lesssim c\frac{\nstoch}{\sqrt{n}}\Vert\bm{\Sigma}\Vert_{F}^2\gamma.
	\end{equation}
	That is
	\begin{displaymath}
	\sqrt{\sum_{j=1}^{\nstoch}\left|\sigma_j^2-\frac{\lambda_j^2}{n}\right|^2}\left/\sum_{j=1}^{\nstoch}\sigma_j^2\right.\lesssim c\frac{\nstoch}{\sqrt{n}}\gamma.
	\end{displaymath}
\end{proof}

\begin{corollary}\label{co:snap_equv}
	Suppose that $\sv_1,\ldots\sv_n$ are random samplings from $\sv$, and the quasimatrix $\widehat{\quasimatrix{U}}^{\samp}_{D\times n}$ is defined as
	\begin{equation}\notag
	\widehat{\quasimatrix{U}}^{\samp}_{D\times n} = \quasimatrix{U}^{\samp}_{D\times n}\diag(\frac{1}{\sqrt{n}},\ldots,\frac{1}{\sqrt{n}}).
	\end{equation}
	Then the number of non-zero singular values of $\widehat{\quasimatrix{U}}^{\samp}_{D\times n}$ is less than or equal to $\nstoch$. Denote the first $\nstoch$ singular values of $\widehat{\quasimatrix{U}}^{\samp}_{D\times n}$ by $\widehat{\lambda}_1\geq\widehat{\lambda}_2\geq\ldots\geq\widehat{\lambda}_{\nstoch}\geq0$, then for large $n$ we have
	\begin{equation}\notag
	\sqrt{\sum_{j=1}^{\nstoch}\left|\sigma_j^2-\widehat{\lambda}_j^2\right|^2}\left/\sum_{j=1}^{\nstoch}\sigma_j^2\right.\lesssim c\frac{\nstoch}{\sqrt{n}}\gamma,
	\end{equation}
	where $\sigma_1\geq\sigma_2\geq\ldots\geq\sigma_{\nstoch}\geq0$ are the singular values of $\quasimatrix{U}^{\trunc}_{D\times \nstoch}$, $c$ is a constant and~$\gamma$ is a random variable.
	
\end{corollary}

\begin{proof}
	Note that $\widehat{\lambda}_j = {\lambda_j}/{\sqrt{n}}$,  the conclusion follows immediately from Theorem~\ref{th:snap_equv}.
\end{proof}

\subsection{A simultaneous low-rank projection approach}
Now let us consider the quasimatrix corresponding to $u^\trunc(\pv,\sv)$. Suppose that the SVD decomposition of $\quasimatrix{U}^{\trunc}_{D\times \nstoch}$ is given by
\begin{equation}\notag
\quasimatrix{U}^{\trunc}_{D\times \nstoch}= {\quasimatrix{U}}{\bm{\Sigma}}{\bm{V}}^*,
\end{equation}
where ${\quasimatrix{U}}$
is a $D\times \nstoch$ quasimatrix with orthonormal columns denoted as $\mathring{u}_1,\ldots, \mathring{u}_{\nstoch}$, ${\bm{\Sigma}}$ is a diagonal matrix with diagonal elements denoted as ${\sigma}_1,\ldots,{\sigma}_{\nstoch}$, and ${\bm{V}}$ is a unitary matrix.
We can then express $u^{\trunc}(\pv,\sv)$ as
\begin{equation}\label{eq:trunc_matrix2}
u^{\trunc}(\pv,\sv) = {\quasimatrix{U}}{\bm{\Sigma}}{\bm{V}}^*(\quasimatrix{S}^{\trunc}_{\Gamma\times\nstoch})^*,
\end{equation}
where
$(\quasimatrix{S}^{\trunc}_{\Gamma\times\nstoch})^*$ is an $\nstoch\times \Gamma$ quasimatrix with orthonormal rows $\gpc_1,\ldots,\gpc_{\nstoch}$. Since ${\bm{V}}$ is a unitary matrix, the rows of the quasimatrix $\bm{\quasimatrix{S}}^*={\bm{V}}^*(\quasimatrix{S}^{\trunc}_{\Gamma\times\nstoch})^*$ are also orthonormal. Denote the rows of $\bm{\quasimatrix{S}}^*$ as $\mathring{\gpc}_1,\ldots,\mathring{\gpc}_{\nstoch}$,  
we can then rewrite~\eqref{eq:trunc_matrix2} as 
\begin{equation}\label{eq:trunc_matrix2bysvd}
u^{\trunc}(\pv,\sv)=  \sum_{j=1}^{\nstoch}{\sigma}_{j}\mathring{u}_j\mathring{\gpc}_j.
\end{equation}
For any $k$, we can consider the partial sum
\begin{equation}\notag
u^{\trunc}_k(\pv,\sv)=  \sum_{j=1}^{k}{\sigma}_{j}\mathring{u}_j\mathring{\gpc}_j.
\end{equation}
It is easy to verify that the relative error of rank $k$ approximation of $u^{\trunc}(\pv,\sv)$ is given by 
\begin{equation}\label{eq:error_partial sum}
\frac{\left\Vert u^{\trunc}(\pv,\sv) - u^{\trunc}_k(\pv,\sv)\right\Vert_{L^2\otimes L^2_{\pdf}}}{\left\Vert u^{\trunc}(\pv,\sv)\right\Vert_{L^2\otimes L^2_{\pdf}}}  = \frac{\left(\sum_{j=k+1}^{\nstoch}{\sigma}_j^2\right)^{1/2}}{\left(\sum_{j=1}^{\nstoch}\sigma_j^2\right)^{1/2}}.
\end{equation}
If we use 
\begin{equation}\notag
X^{(k)} = \mathrm{span}\{\mathring{u}_1,\ldots,\mathring{u}_k\},\ S^{(k)} =\mathrm{span} \{\mathring{\gpc}_1,\ldots,\mathring{\gpc}_k\}
\end{equation}
as the low-rank approximation subsapces of the physical and stochastic spaces respectively, we may result in a much smaller linear system compared to~\eqref{eq:linsg}. However, since the $D\times \nstoch$ quasimatrix $\quasimatrix{U}^{\trunc}_{D\times \nstoch}$ cannot be determined in advance, some adjustments are necessary to adapt this idea to design a solver for~\eqref{eq:spde}.

It is worth noting that the accuracy of the rank $k$ approximation of~$u^{\trunc}(\pv,\sv)$  depends only on the quasimatrix $\quasimatrix{U}^{\trunc}_{D\times \nstoch}$. In practical implementations,  quadrature formulas are typically used to compute the integral in the Gram–Schmidt orthogonalization process for the QR factorization of a quasimatrix. This inspires us to construct the low-rank basis functions of the stochastic space (i.e., basis functions of $S^{(k)}$) and to determine the rank $k$ as follows: 

First, we compute a numerical solution on a coarse grid of the physical space. If the coefficients of the gPC basis functions are denoted by $u_1^{\coarse}(\pv),\ldots, u_{\nstoch}^{\coarse}(\pv)$, the solution on the coarse grid of the physical space can be expressed as:
\begin{equation}\notag
u^{\coarse}(\pv,\sv) = \quasimatrix{U}^{\coarse}_{D\times\nstoch}(\quasimatrix{S}^{\trunc}_{\Gamma\times \nstoch})^*, 
\end{equation}
where
\begin{equation}\notag
\quasimatrix{U}^{\coarse}_{D\times\nstoch} = \left[\begin{array}{llll}
\Big|    	&   \Big|  	&        	    &  \Big|        \\ 
u_1^{\coarse} 		& u_2^{\coarse} 		& \ldots 		& u_{\nstoch}^{\coarse} \\
\Big|    	&   \Big|  	&        	    &  \Big|      
\end{array}\right].  
\end{equation}
Next, we perform the SVD decomposition of $\quasimatrix{U}^{\coarse}_{D\times\nstoch}$:
\begin{equation}\label{eq:coarse_svd}
\quasimatrix{U}^{\coarse}_{D\times\nstoch} = \widehat{\quasimatrix{U}}\widehat{\bm{\Sigma}}\widehat{\bm{V}}^*,
\end{equation}
where $\widehat{\quasimatrix{U}}$ is a $D\times\nstoch$ quasimatrix and $\widehat{\bm{\Sigma}}$ is a diagonal matrix. The diagonal elements of $\widehat{\bm{\Sigma}}$, denoted by $\widehat{\sigma}_1,\ldots,\widehat{\sigma}_{\nstoch}$, approximate $\sigma_1,\ldots,\sigma_{\nstoch}$ form~\eqref{eq:trunc_matrix2bysvd}. Defining 
\begin{equation}\label{eq:stoch_basis_quasi}
\widehat{\quasimatrix{S}}^{*} =\widehat{\bm{V}}^*(\quasimatrix{S}^{\trunc}_{\Gamma\times\nstoch})^*,
\end{equation}
we note that $\widehat{\quasimatrix{S}}^{*}$ is a $\Gamma\times\nstoch$ quasimatrix with orthonormal rows. Denoting the rows of $\widehat{\quasimatrix{S}}^{*}$ by $\widehat{\gpc}_1(\sv),\ldots,\widehat{\gpc}_{\nstoch}(\sv)$, the low-rank approximation subspace of stochastic space is given by:  
\begin{equation}\notag
S^{(k)} = \mathrm{span}\{\widehat{\gpc}_1,\ldots,\widehat{\gpc}_k\},
\end{equation}
where $k$ is the smallest number such that 
\begin{equation}\label{eq:find_k}
\left(\sum_{j=k+1}^{\nstoch}(\widehat{\sigma}_j)^2\right)^{1/2}\left/ \left(\sum_{j=1}^{\nstoch}(\widehat{\sigma}_j)^2\right)^{1/2}\right. \leq \tol.
\end{equation}
Here, $\tol$ is the desired accuracy.

For the low-rank basis functions of the physical space (i.e., basis functions of $X^{(k)}$), choosing the columns of the quasimatrix~$\widehat{\quasimatrix{U}}$ becomes less effective when the mesh grid is refined. Instead, we construct the low-rank basis functions of the physical space based on Theorem~\ref{th:snap_equv}. This theorem asserts that if the number of snapshots (i.e., $u^{\trunc}(\pv,\sv_i),\ i=1,\ldots,n$) is sufficiently large, the non-zero singular values of ${\quasimatrix{U}}^{\samp}_{D\times n}$ can 
be arbitrarily close to the singular values of $\quasimatrix{U}^{\trunc}_{D\times \nstoch}$ with a scale factor. In other words, the range of~$\quasimatrix{U}^{\samp}_{D\times n}$ aligns with the range of~$\quasimatrix{U}^{\trunc}_{D\times \nstoch}$. This suggests that the left singular vectors of~$\quasimatrix{U}^{\samp}_{D\times n}$  can be used in place of the left singular vectors of $\quasimatrix{U}^{\trunc}_{D\times \nstoch}$ as the basis functions of the physical space.

To further enhance efficiency, we employ a greedy algorithm with a residual-free indicator~\cite{Chen2019resfree} to identify the most important $k$ sampling points from a set of random samples $\bm{\Theta}$. For computational efficiency, both snapshots and residual-free indicators are computed on a coarse grid. The identified $k$ sampling points are denoted as~$\{\sv_1^{\coarse},\ldots\sv_k^{\coarse}\}$, and we then compute $u(\pv,\sv_i^{\coarse})$ for $i=1,\ldots,k$ on the fine grid. These fine grid snapshots are denoted as $u^{\fine}(\pv,\sv_i^{\coarse})$. Suppose that the SVD decomposition of the quasimatrix $\quasimatrix{U}^{\fine}$ is given by
\begin{equation}\notag
\quasimatrix{U}^{\fine} = \widetilde{\quasimatrix{U}}\widetilde{\bm{\Sigma}}\widetilde{\bm{V}}^*,
\end{equation}
where
\begin{equation}\notag
\quasimatrix{U}^{\fine} = \left[\begin{array}{llll}
\Big|    	&   \Big|  	&        	    &  \Big|        \\ 
u^{\fine}(\pv,\sv_1^{\coarse}) 		& u^{\fine}(\pv,\sv_2^{\coarse}) 		& \ldots 		& u^{\fine}(\pv,\sv_k^{\coarse}) \\
\Big|    	&   \Big|  	&        	    &  \Big|      
\end{array}\right].
\end{equation}
We then use 
\begin{equation}\notag
X ^{(k)} = \mathrm{span}\{\widetilde{u}_1, \widetilde{u}_2 , \ldots , \widetilde{u}_k \}
\end{equation}
as the low-rank approximation subspace of the physical space, where $k$ is the same number as in~\eqref{eq:find_k}, and $\widetilde{u}_1, \widetilde{u}_2 , \ldots , \widetilde{u}_k$ are the columns of $\widetilde{\quasimatrix{U}}$.

Suppose that $\bm{W}^{(k)}$ and $\bm{V}^{(k)}$ are matrices satisfying
\begin{equation}\label{eq:coordinates}
\widetilde{\quasimatrix{U}} =  \quasimatrix{X}\bm{W}^{(k)}, \ \widetilde{\quasimatrix{S}}^*   =(\bm{V}^{(k)})^*(\quasimatrix{S}^{\trunc}_{\Gamma\times\nstoch})^*,
\end{equation}
where
\begin{equation}\notag
\quasimatrix{X} = \left[\begin{array}{llll}
\Big|    	&   \Big|  	&        	    &  \Big|        \\ 
\pbf_1(\pv) 		& \pbf_2(\pv) 		& \ldots 		& \pbf_{\nphy}(\pv) \\
\Big|    	&   \Big|  	&        	    &  \Big|      
\end{array}\right], \ \widetilde{\quasimatrix{S}}^* = \left[\begin{array}{ccc}
\text{---} 	    &  \widehat{\gpc}_1	&  \text{---}           \\ 
\text{---}      &  \widehat{\gpc}_2 		&  \text{---}\\
&   \vdots  	&       \\      
\text{---}  	&  \widehat{\gpc}_{k} & \text{---}
\end{array}\right],
\end{equation}
and the matrix $\bm{V}^{(k)}$ consists of the first $k$ columns of $\widehat{\bm{V}}$ defined in~\eqref{eq:coarse_svd}. Note that the columns of the quasimatrix $\quasimatrix{X}$ correspond to the physical basis functions. Consequently, the columns of $\bm{W}^{(k)}$ represent the coordinates of $\widetilde{u}_1,\widetilde{u}_2,\ldots \widetilde{u}_k$ with respect to the physical basis functions $\{\pbf_s(\pv)\}_{s=1}^{\nphy}$. Similarly, the rows of the quasimatrix $(\quasimatrix{S}^{\trunc}_{\Gamma\times\nstoch})^*$ correspond to the stochastic basis functions, and the rows of $(\bm{V}^{(k)})^*$ represent the coordinates of $\widehat{\gpc}_1, \widehat{\gpc}_2, \ldots, \widehat{\gpc}_k$ with respect to the stochastic basis functions~$\{\gpc_s(\pv)\}_{s=1}^{\nstoch}$.

If we use 
\begin{equation}\notag
   W^{(k)} =  X^{(k)}\otimes  S^{(k)}:= \mathrm{span} \left\{\widetilde{u}_s(\pv)\widehat{\gpc}_j(\sv)\left| \widetilde{u}_s\in X^{(k)}, \widehat{\gpc}_j\in S^{(k)} \right.\right\}
\end{equation}
as the approximation subspace, where
\begin{equation}\label{eq:phy_stoch_lowrank}
X ^{(k)} = \mathrm{span}\{\widetilde{u}_1, \widetilde{u}_2 , \ldots , \widetilde{u}_k \},\ S^{(k)} = \mathrm{span}\{\widehat{\gpc}_1,\ldots,\widehat{\gpc}_k\},
\end{equation}
we obtain a significantly smaller linear system compared to the full system~\eqref{eq:linsg}. We refer to this reduced system as the {\itshape low-rank projected linear system}, as it results from projecting the original system~\eqref{eq:linsg} onto the low-rank subspace~$W^{(k)}$.

By substituting the low-rank projected approximation of the solution, i.e., 
\begin{equation}\label{eq:lowrank_app}
    u^{\trunc}_{k}(\pv,\sv): =\sum_{j=1}^{k}\sum_{s=1}^{k}x_{sj}^{(k)}\widetilde{u}_s(\pv)\widehat{\gpc}_j(\sv),
\end{equation}
into the weak formulation~\eqref{eq:weak}, and using the coordinate representation~\eqref{eq:coordinates}, we obtain the low-rank projected linear system in matrix form:
\begin{equation}\label{eq:SimLRP}
    \sum_{i=1}^{K} \bm{A}_i^{(k)} \bm{X}^{(k)}(\bm{G}_i^{(k)})^\trans = \bm{f}^{(k)} (\bm{h}^{(k)})^\trans,
\end{equation}
where $\bm{X}^{(k)}$ collects the unknown coefficients $x_{sj}^{(k)}$ from~\eqref{eq:lowrank_app}, with entries defined by
\begin{equation}\notag
    \bm{X}^{(k)}(s,j)=x_{sj}^{(k)}.
\end{equation}
The reduced matrices and vectors in~\eqref{eq:SimLRP} are given by:
\begin{alignat}{2}
    \bm{A}_i^{(k)}&=(\bm{W}^{(k)})^*\bm{A}_i\bm{W}^{(k)}, \ &\bm{f}^{(k)}&=(\bm{W}^{(k)})^*\bm{f},\label{eq:update_phy}\\
    \bm{G}_i^{(k)}&= (\bm{V}^{(k)})^\trans\bm{G}_i\left((\bm{V}^{(k)})^*\right)^{\trans}, \ &\bm{h}^{(k)}&=(\bm{V}^{(k)})^{\trans}\bm{h},\label{eq:update_stoch}
\end{alignat}
where the matrices and vectors on the right-hand sides of~\eqref{eq:update_phy}--\eqref{eq:update_stoch} are defined in~\eqref{eq:stoch_mat}--\eqref{eq:phy_mat} and~\eqref{eq:coordinates}.

Note that the low-rank projected matrix equation~\eqref{eq:SimLRP} involves only $k\times k$ unknowns, and can therefore be solved much more efficiently than the full matrix equation~\eqref{eq:linsg_matrix}. Once  $\bm{X}^{(k)}$  is solved from~\eqref{eq:SimLRP}, the low-rank approximation of the solution can be reconstructed as
\begin{equation}\label{eq:construct_lowrank}
    u^{\trunc}_{k}(\pv,\sv) = \widetilde{\quasimatrix{U}}\bm{X}^{(k)}\widetilde{\quasimatrix{S}}^* =  \quasimatrix{X}\bm{W}^{(k)}\bm{X}^{(k)}(\bm{V}^{(k)})^*(\quasimatrix{S}^{\trunc}_{\Gamma\times\nstoch})^*.
\end{equation}
On the other hand, the full solution obtained from the standard stochastic Galerkin method, as given in equation~\eqref{eq:appsg}, can be written in the form
\begin{equation}\label{eq:construct_full}
    u^{\ap}(\pv,\sv) = \quasimatrix{X}\bm{U}(\quasimatrix{S}^{\trunc}_{\Gamma\times\nstoch})^*.
\end{equation}
Since both~\eqref{eq:construct_lowrank} and~\eqref{eq:construct_full} approximate the solution to problem~\eqref{eq:spde}, the matrix product $\bm{W}^{(k)}\bm{X}^{(k)}(\bm{V}^{(k)})^*$ can be viewed as a low-rank approximation to the full coefficient matrix $\bm{U}$.

In practical implementations, all operations are conducted in matrix form. In this study, the meshgrid in the physical space is uniform, and the $L^2$ norm of a function is approximated by the~$2$-norm of the corresponding vector multiplied by a factor of~$1/\nphy$. Under these conditions, we employ matrix SVD to replace quasimatrix SVD decomposition. Algorithm~\ref{alg:SimLRP} provides the pseudocode for the SimLRP approach.

\begin{algorithm}[h]
	\caption{A simultaneous low-rank projection (SimLRP) approach}\label{alg:SimLRP}
	\begin{algorithmic}[1]
		\State
		\textbf{INPUT:}  $p$: the gPC oder; $\tol$: the tolerance in\eqref{eq:find_k}; $T^{\coarse}$: a coarse physical grid; $T$: a fine physical grid; $\bm{\Theta}$:~a set of random samples.
		\State	Generate  $\bm{A}_i^{\coarse},\bm{f}^{\coarse}$ for the coarse physical gird  $T^{\coarse}$, and  $\bm{A}_i,\bm{f}$  for the fine physical gird  $T$, generate $\bm{G}_i$ and $\bm{h}$.
		\State Solve the matrix equation
		\begin{equation}\notag
		\sum_{i=1}^{K} \bm{A}_i^{\coarse} \bm{U}^{\coarse}\bm{G}_i^\trans = \bm{f}^{\coarse} \bm{h}^\trans,
		\end{equation}
		\State Compute SVD decomposition of $\bm{U}^{\coarse}$, i.e., $\bm{U}^{\coarse}=\widehat{\bm{U}}\widehat{\bm{\Sigma}}\widehat{\bm{V}}^*$.
		\State $\ k$ $\longleftarrow$ the smallest number satisfying~\eqref{eq:find_k}, where $\widehat{\sigma}_1,\ldots,\widehat{\sigma}_{\nstoch}$ are the diagonal of $\widehat{\bm{\Sigma}}$.
		
		\State $\bm{V}^{(k)}$ $\longleftarrow$ the first $k$ columns of $\widehat{\bm{V}}$.
		
		\State  $\ \bm{\theta}$ $\longleftarrow$ the most important $k$ samples selected from $\bm{\Theta}$.

		\State $\bm{U}^{\fine}$ $\longleftarrow$ snapshots on the fine physical grid corresponding to $\bm{\theta}$.
		\State $\bm{W}^{(k)}$ $\longleftarrow$ the first $k$ left singular vectors of $\bm{U}^{\fine}$.
		
		\State $\bm{X}^{(k)}$$\longleftarrow$ solve~\eqref{eq:SimLRP}.
		\State Construct low-rank approximation by~\eqref{eq:construct_lowrank}. 
	\end{algorithmic}
\end{algorithm}

\subsection{An alternating low-rank projection approach}
While the SimLRP approach provides the capability to simultaneously reduce the dimensionality of both the physical and stochastic spaces, the tolerance indicator~\eqref{eq:find_k} may deteriorate for smaller values of $\tol$. Drawing inspiration from~\cite{jain2013low,lee2022enhanced}, we propose an alternating low-rank projection (AltLRP) approach to mitigate the deterioration of tolerance indicators for smaller values.

The idea is quite straightforward. After applying SimLRP, we alternatively implement the low-rank projection in the stochastic and physical spaces to refine and update the low-rank approximation. Suppose we have obtained $\bm{W}^{(k)}$, $\bm{X}^{(k)}$, $\bm{V}^{(k)}$ from the SimLRP approach. We first update $\bm{W}^{(k)}$ and $\bm{V}^{(k)}$ based on the approximate solution obtained from~\eqref{eq:construct_lowrank} in the SimLRP approach. We then use
    \begin{equation}\label{eq:low_rank_subspace_all}
       W_{\mathrm{S}}^{(k)}=X_h\otimes  S^{(k)}:= \mathrm{span} \left\{\pbf_s(\pv)\widehat{\gpc}_j(\sv)\left| \pbf_s\in X_h, \widehat{\gpc}_j\in S^{(k)} \right.\right\}
    \end{equation}
as the low-rank approximation subspace, where $X_h$ and $S^{(k)}$ are defined in~\eqref{eq:phy_stoch_full} and~\eqref{eq:phy_stoch_lowrank} respectively. Similar to the SimLRP approach, this leads to a low-rank projected linear system in matrix form:
\begin{equation}\label{eq:AltLRP_stoch}
    \sum_{i=1}^{K} \bm{A}_i \bm{Y}^{(k)}(\bm{G}_i^{(k)})^\trans = \bm{f} (\bm{h}^{(k)})^\trans,\\
\end{equation}
where $\bm{Y}^{(k)}$ collects the unknown coefficients, and the other matrices and vectors are defined in~\eqref{eq:phy_mat} and~\eqref{eq:update_stoch}.

Once $\bm{Y}^{(k)}$ is solved from~\eqref{eq:AltLRP_stoch}, the low-rank approximation solution, similar to the SimLRP approach, is given by
\begin{equation}\label{eq:construct_partial_stoch}
    \quasimatrix{X}\bm{Y}^{(k)}\widetilde{\quasimatrix{S}}^* = \quasimatrix{X}\bm{Y}^{(k)}(\bm{V}^{(k)})^*(\quasimatrix{S}^{\trunc}_{\Gamma\times\nstoch})^*.
\end{equation}
We then update $\bm{W}^{(k)}$ and $\bm{V}^{(k)}$ based on the approximation solution obtained from~\eqref{eq:construct_partial_stoch}, and use
 \begin{equation}\notag
    W_{\mathrm{P}}^{(k)}=X^{(k)}\otimes  S_p:= \mathrm{span} \left\{\widetilde{\pbf}_s(\pv){\gpc}_j(\sv)\left| \widetilde{\pbf}_s\in X^{(k)}, {\gpc}_j\in S_p \right.\right\}
\end{equation}
as the low-rank approximation subspace,  where $X^{(k)}$ and $S_p$ are defined in~\eqref{eq:phy_stoch_lowrank} and~\eqref{eq:phy_stoch_full}, respectively. Similar to the SimLRP approach, this leads to a low-rank projected linear system in matrix form:
\begin{equation}\label{eq:AltLRP_phy}
    \sum_{i=1}^{K} \bm{A}_i^{(k)} \bm{Z}^{(k)}\bm{G}_i^\trans = \bm{f}^{(k)} \bm{h}^\trans. 
\end{equation}
where $\bm{Z}^{(k)}$ collects the unknown coefficients, and the other matrices and vectors are defined in~\eqref{eq:stoch_mat} and~\eqref{eq:update_phy}.

Once $\bm{Z}^{(k)}$ is solved from~\eqref{eq:AltLRP_phy}, the low-rank approximation solution, similar to the SimLRP approach, is given by
\begin{equation}\label{eq:construct_partial_phy}
    \widetilde{\quasimatrix{U}}\bm{Z}^{(k)}(\quasimatrix{S}^{\trunc}_{\Gamma\times\nstoch})^* = \quasimatrix{X}\bm{W}^{(k)}\bm{Z}^{(k)}(\quasimatrix{S}^{\trunc}_{\Gamma\times\nstoch})^*.
\end{equation}
We can then update $\bm{W}^{(k)}$ and $\bm{V}^{(k)}$ based on the approximate solution obtained from~\eqref{eq:construct_partial_phy}, and use $W_{\mathrm{S}}^{(k)}$ defined in \eqref{eq:low_rank_subspace_all} as the low-rank approximation subspace.

By repeating this procedure, we obtain the alternating low-rank projection (AltLRP) approach. The method is named as such because it alternates between applying low-rank projections in the stochastic and physical spaces. To reduce the number of iterations in solving~\eqref{eq:AltLRP_stoch} and~\eqref{eq:AltLRP_phy}, the solution from the previous step can be used as an initial guess.

\begin{algorithm}[h]
	\caption{Alternating low-rank projection (AltLRP) approach}\label{alg:AltLRP}
	\begin{algorithmic}[1]
		
		\State Implement Algorithm~\ref{alg:SimLRP}, and keep the notations.
		
		\State  [$\bm{W}^{(k)}$,$\bm{\Sigma}^{(k)}$,$\bm{V}^{(k)}$] $\longleftarrow$\Call{TruncatedSVD}{$\bm{X}^{(k)}$, $\bm{W}^{(k)}$, $\bm{V}^{(k)}$}.

		\For {$i=1,2,\ldots i_{\max}$}

		\State Update $\bm{G}_i^{(k)}$ and $\bm{h}^{(k)}$ by~\eqref{eq:update_stoch}.
		
		\State $\bm{Y}^{(k)}$ $\longleftarrow$ solve~\eqref{eq:AltLRP_stoch} with initial guess $\bm{W}^{(k)}\bm{\Sigma}^{(k)}$. 
		
		\State [$\bm{W}^{(k)}$,$\bm{\Sigma}^{(k)}$,$\bm{V}^{(k)}$] $\longleftarrow$\Call{TruncatedSVD}{$\bm{Y}^{(k)}$, $\bm{E}$, $\bm{V}^{(k)}$}.

		\State Update $\bm{A}_i^{(k)}$ and $\bm{f}^{(k)}$ by~\eqref{eq:update_phy}.
		
		\State $\bm{Z}^{(k)}$ $\longleftarrow$ solve~\eqref{eq:AltLRP_phy} with initial guess $\bm{\Sigma}^{(k)}(\bm{V}^{(k)})^*$. 
		
		\State [$\bm{W}^{(k)}$,$\bm{\Sigma}^{(k)}$,$\bm{V}^{(k)}$] $\longleftarrow$\Call{TruncatedSVD}{$\bm{Z}^{(k)}$, $\bm{W}^{(k)}$, $\bm{E}$}.
		\EndFor
		\State Update $\bm{A}_i^{(k)}$, $\bm{f}^{(k)}$, $\bm{G}_i^{(k)}$,  $\bm{h}^{(k)}$ by~\eqref{eq:update_phy}--\eqref{eq:update_stoch}.
		\State $\bm{X}^{(k)}$$\longleftarrow$ solve~\eqref{eq:SimLRP} with initial guess $\bm{\Sigma}^{(k)}$.
		
		\State Construct low-rank approximation by~\eqref{eq:construct_lowrank}. 
	\end{algorithmic}
\end{algorithm}

\begin{algorithm}[!h]
	\caption{Truncated SVD decomposition}  
	\begin{algorithmic}[1]
		\Function{[$\bm{W}^{(k)}$,$\bm{\Sigma}^{(k)}$, $\bm{V}^{(k)}$] = TruncatedSVD}{$\bm{X}$,$\bm{W}^{(k)}$, $\bm{V}^{(k)}$}
		
		\State Compute the SVD decomposition of $\bm{X}$, i.e., $\bm{X}=\bm{W}\bm{\Sigma}\bm{V}^*$.
		\State $\bm{W}^{(k)}_{\temp}$ $\longleftarrow$ first  $k$ columns of $\bm{W}$.
		
		\State $\bm{V}^{(k)}_{\temp}$ $\longleftarrow$ first  $k$ rows of $\bm{V}$.
		
		\State $\bm{\Sigma}^{(k)}_{\temp}$ $\longleftarrow$ first  $k\times k$ leading principal submatrix of $\bm{\Sigma}$.
		
		\State $\bm{W}^{(k)}$$\longleftarrow$ $\bm{W}^{(k)}\bm{W}^{(k)}_{\temp}$, $\bm{V}^{(k)}$$\longleftarrow$ $\bm{V}^{(k)}\bm{V}^{(k)}_{\temp}$, $\bm{\Sigma}^{(k)}$$\longleftarrow$ $\bm{\Sigma}^{(k)}_{\temp}$.

		\EndFunction
	\end{algorithmic}
	\label{alg:stage_aem}
\end{algorithm}

Algorithm~\ref{alg:AltLRP} provides the pseudocode for the AltLRP approach. Since $\bm{W}^{(k)}$, $\bm{X}^{(k)}$, $\bm{V}^{(k)}$ obtained from the SimLRP approach are already quite accurate, typically only one iteration (i.e., $i_{\max} = 1$) is required in the AltLRP approach to achieve the desired accuracy.

\section{Numerical results}\label{sec:numericaltests}

In this section, we consider two problems: a diffusion problem and a Helmholtz problem. For both test problems, spatial discretization is performed using $Q_1$ finite elements based on the IFISS package~\cite{ifiss-siamreview}. All the results presented in this study are obtained using MATLAB R2016a on a desktop computer with 2.90GHz Intel Core i7-10700 CPU.

To assess accuracy, we define the relative error as follows:
\begin{equation}\label{eq:defRelErr}
\mathtt{RelErr} = \frac{\left\Vert\bm{W}^{(k)}\bm{X}^{(k)}\bm{V}^{(k)} - \bm{U}\right\Vert_{F}}{\left\Vert \bm{U}\right\Vert_{F}},
\end{equation}
where $\bm{W}^{(k)}$, $\bm{X}^{(k)}$, $\bm{V}^{(k)}$ are the matrices corresponding to the low-rank approximation with rank~$k$, and $\bm{U}$ is the matrix in~\eqref{eq:linsg_matrix}.

\subsection{Test problem 1}
In this problem, we investigate the diffusion equation with random inputs, given by
\begin{equation}\notag
\left\{
\begin{aligned}
-\nabla\cdot(a(\pv,\sv)\nabla u(\pv,\sv)) = 1 &\quad \mbox{in}\quad D\times\Gamma,\\
u(\pv,\sv) = 0 &\quad \mbox{on}\quad \partial D\times \Gamma,\\
\end{aligned}\right.
\end{equation}
where $D = [-1,1]^2$ is the spatial domain, and $\partial D$ represents the boundary of $D$. The diffusion coefficient $a(\pv,\sv)$ is modeled as a truncated Karhunen--Lo\`eve (KL) expansion~\cite{Ghanem2003,Elman2007} of a random field with a mean function $a_0(\pv)=0.2$, a standard deviation $\sigma=0.1$, and the covariance function $\mathrm{Cov}\,(\bm{x},\bm{y})$ given by
\begin{equation}\notag
\mathrm{Cov}\,(\bm{x},\bm{y})=\sigma^2 \exp\left(-\frac{|x_1-y_1|}{c}-\frac{|x_2-y_2|}{c}
\right),\label{covariance}
\end{equation}
where $\pv=[x_1,x_2]^T$, $\bm{y}=[y_1,y_2]^T$ and $c=1$ is the correlation length.
The KL expansion takes the form
\begin{equation}\label{eq:kl}\notag
a(\pv,\sv)=a_0(\pv)+\sum_{i=1}^{N}a_i(\pv)\svc_i=a_0(\pv)+\sum_{i=1}^{N}\sqrt{\lambda_i}c_i(\pv)\svc_i,
\end{equation}
where $\{\lambda_i,c_i(\pv)\}_{i=1}^{N}$ are the eigenpairs of  $\mathrm{Cov}\,(\pv,\bm{y})$,
$\{\svc_i\}^{N}_{i=1}$ are uncorrelated random variables, and
$N$ is the number of KL modes retained. For this test problem, we assume that the random variables  $\{\svc_i\}^{N}_{i=1}$ are identically independent distributed uniform random variables on $[-1, 1]$.

In this test problem, the gPC order $p$ is set to $4$, the degrees of freedom (DOF) of the coarse grid in physical space (i.e., $\bm{T}^c$) is set to $\nphy = 33\times 33$, and the cardinality of the random set~$\bm{\Theta}$ is set to $|\bm{\Theta}| = \min(2\nstoch,\nphy)$. All the linear systems, i.e.,~\eqref{eq:SimLRP} and~\eqref{eq:AltLRP_stoch}--\eqref{eq:AltLRP_phy}, are solved using the preconditioned conjugate gradients (PCG) method  with a mean-based preconditioner~\cite{Powell2009} and a tolerance of $10^{-8}$.

\begin{figure}[htbp]
	\begin{center}
		\subfloat[$N=5$]{
			\includegraphics[width=0.27\linewidth]{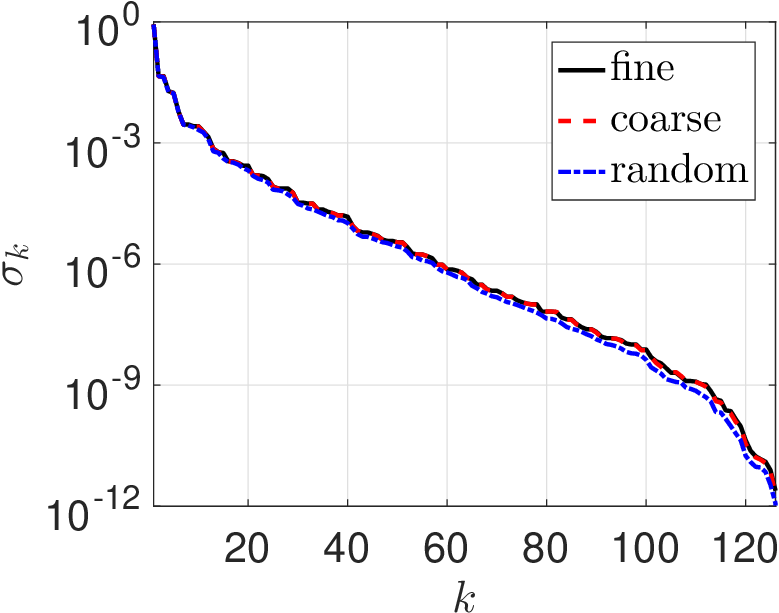}
		}\qquad
		\subfloat[$N=7$]{
			\includegraphics[width=0.27\linewidth]{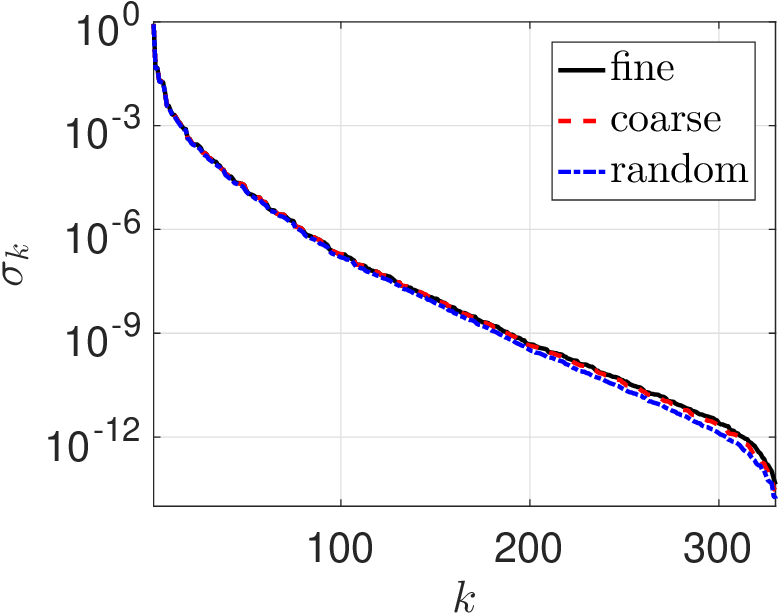}
		}\qquad
		\subfloat[$N=10$]{
			\includegraphics[width=0.27\linewidth]{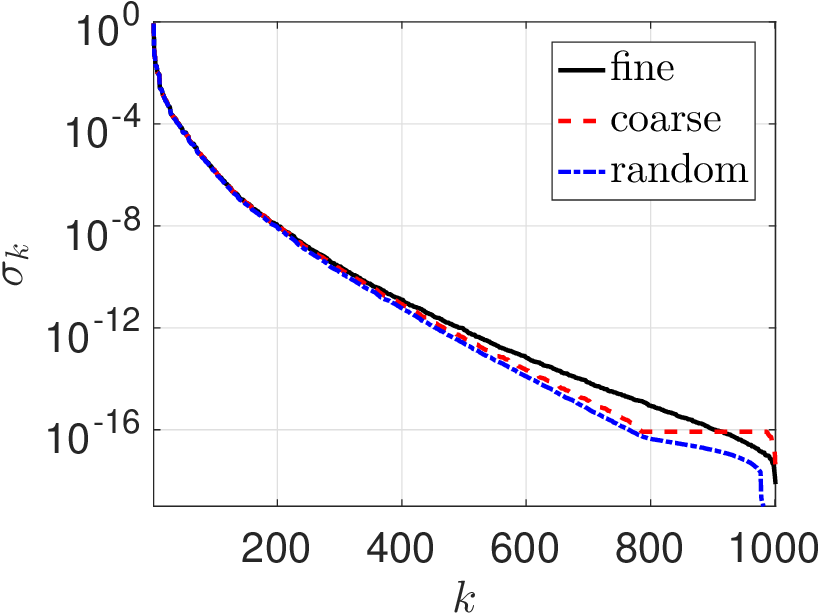}
		}
	\end{center}

	\caption{Singular values for different $N$.}\label{fig:test1SVD}
	{\small In this figure, `fine' refers to the singular values of $\quasimatrix{U}^{\trunc}_{D\times \nstoch}$ with $\nphy= 257^2$; `coarse' refers to the singular values of $\quasimatrix{U}^{\coarse}_{D\times \nstoch}$ with $\nphy= 33^2$;  and `random' refers to the singular values of $\widehat{\quasimatrix{U}}^{\samp}_{D\times n}$ in Corollary~\ref{co:snap_equv} with $|\bm{\Theta}|=\min(2\nstoch,\nphy)$ and $\nphy=33^2$.}
\end{figure}

In Figure~\ref{fig:test1SVD}, we plot the singular values for different $N$. The DOF of the physical space is $\nphy = 33^2$ for the coarse grid and $\nphy = 257^2$ for the fine grid. It is noteworthy that, in this example, the $L^2$ norm of a function is simply approximated by  the $2$-norm of its corresponding vector, multiplied by a factor of $1/\nphy$, since the mesh grid in physical space is uniform.  From the figure, we observe that the singular values for `fine', `coarse' and `random' coincide very well with each other, reflecting the conclusion of Theorem~\ref{th:snap_equv}. It also demonstrates that the singular values of~$\quasimatrix{U}^{\trunc}_{D\times \nstoch}$ are not very sensitive to the DOF of the physical mesh grid, and thus allowing the use of singular values computed on a coarse grid to determine the rank $k$ for a desired accuracy $\tol$.

\begin{figure}[htbp]
	\begin{center}
		\subfloat[$N=5,\ \nphy = 129^2$]{
			\includegraphics[width=0.27\linewidth]{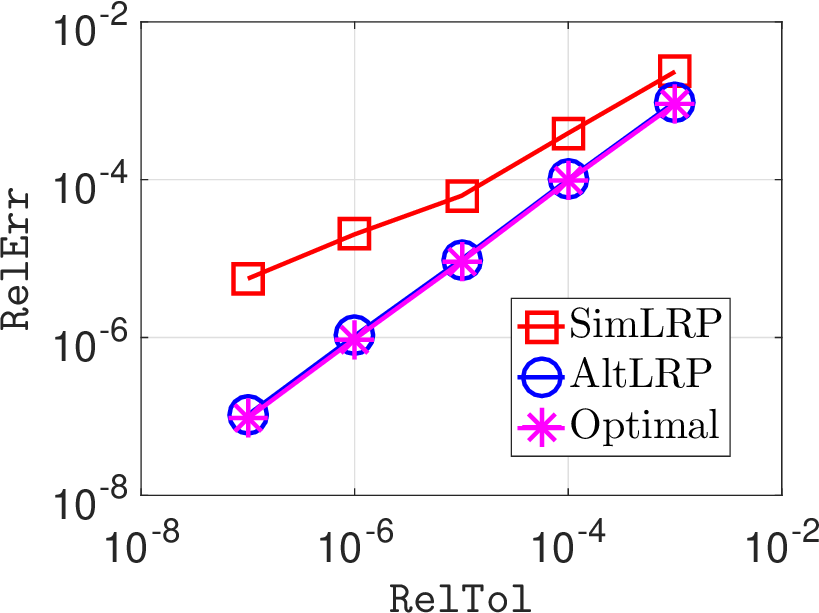}
		}\qquad
		\subfloat[$N=7,\ \nphy = 129^2$]{
			\includegraphics[width=0.27\linewidth]{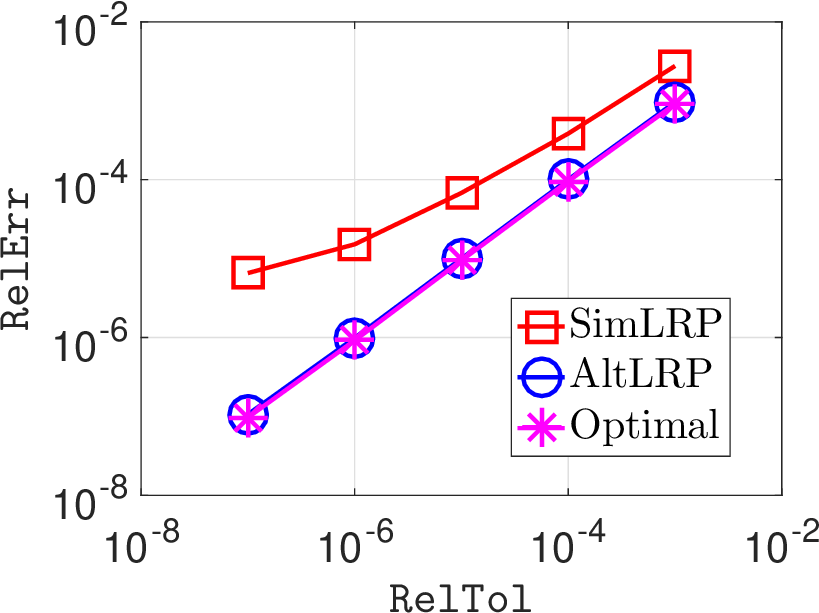}
		}\qquad
		\subfloat[$N=10,\ \nphy = 129^2$]{
			\includegraphics[width=0.27\linewidth]{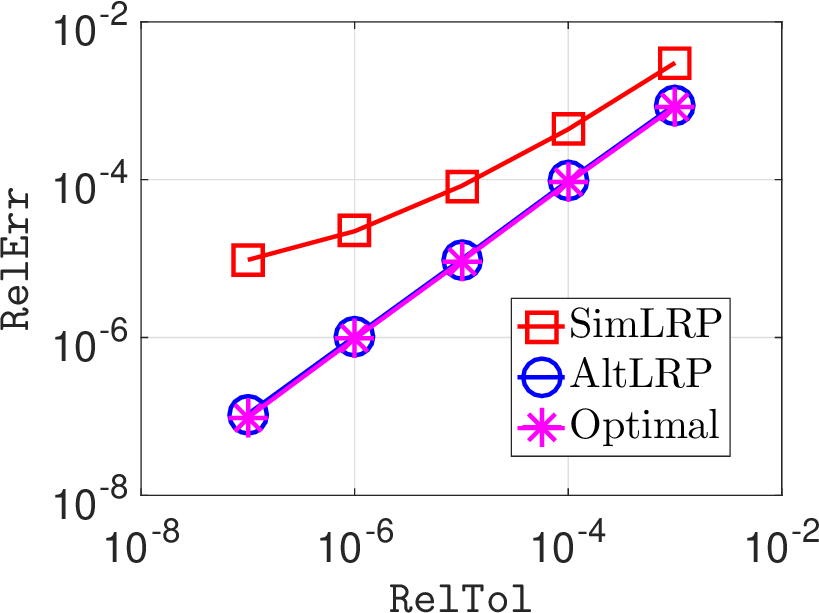}
		}\\
		\subfloat[$N=5, \ \nphy = 257^2$]{
			\includegraphics[width=0.27\linewidth]{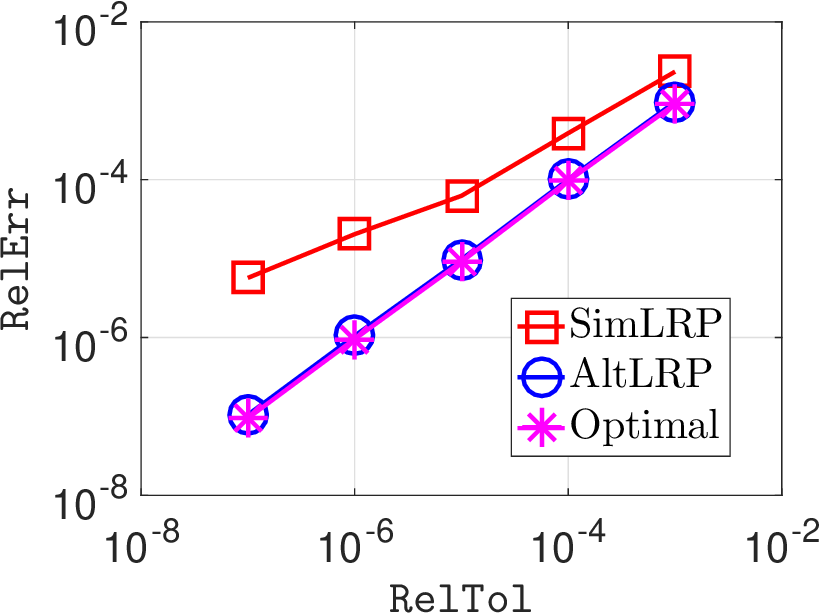}
		}\qquad
		\subfloat[$N=7,\ \nphy = 257^2$]{
			\includegraphics[width=0.27\linewidth]{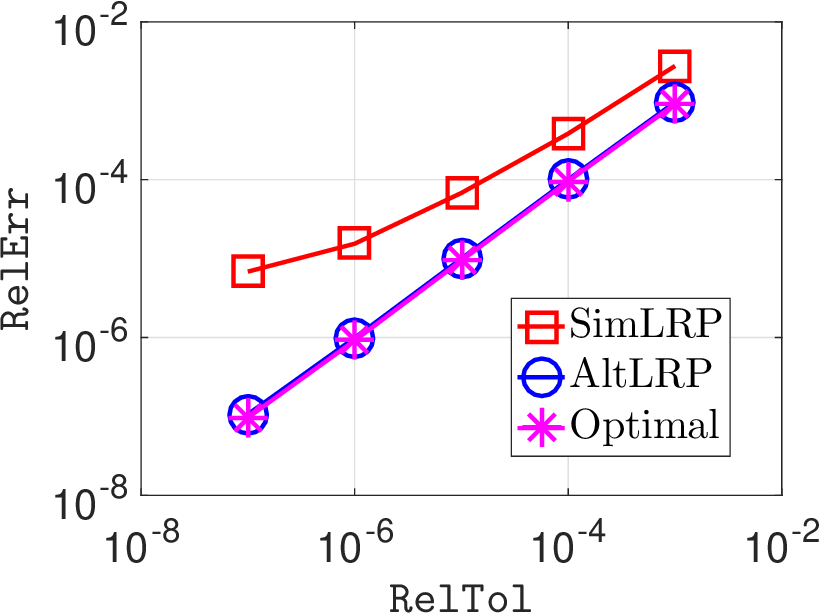}
		}\qquad
		\subfloat[$N=10,\ \nphy = 257^2$]{
			\includegraphics[width=0.27\linewidth]{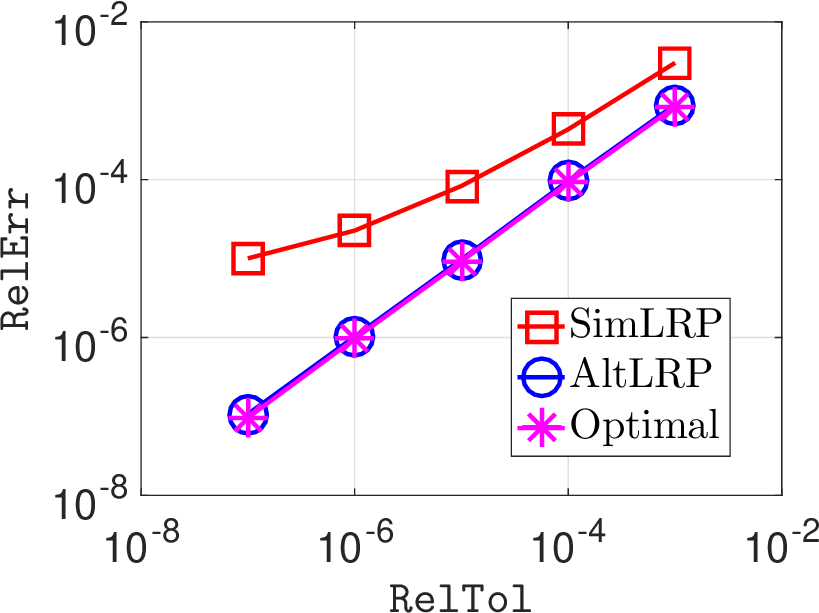}
		}
	\end{center}
	\caption{Relative errors with respect to different $\tol$.}\label{fig:test1b}
\end{figure}

In Figure~\ref{fig:test1b}, we present the relative errors concerning different values of $\tol$, where the relative error is defined by~\eqref{eq:defRelErr}. The `optimal' low-rank approximation of the solution is obtained through the SVD decomposition of $\bm{U}$ with the corresponding $k$ terms retained. In the AltLRP approach, the value of $i_{\max}$ is set to be~$1$. From the figure, it is evident that the AltLRP approach performs nearly optimally. While the SimLRP approach also demonstrates good performance, the tolerance indicator~\eqref{eq:find_k} for it deteriorates for smaller values of $\tol$.

\begin{figure}[htbp]
	\begin{center}
		\subfloat[$N=5,\ \nphy = 129^2$]{
			\includegraphics[width=0.27\linewidth]{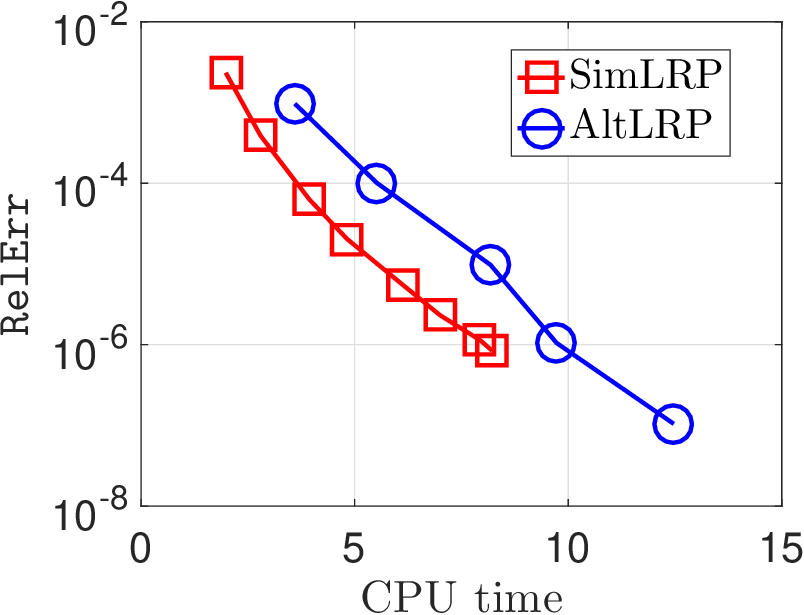}
		}\qquad
		\subfloat[$N=7,\ \nphy = 129^2$]{
			\includegraphics[width=0.27\linewidth]{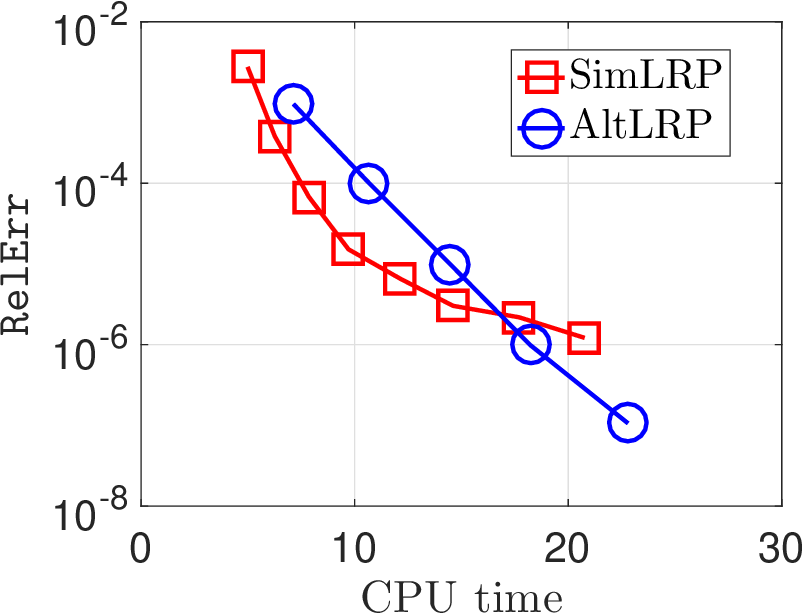}
		}\qquad
		\subfloat[$N=10,\ \nphy = 129^2$]{
			\includegraphics[width=0.27\linewidth]{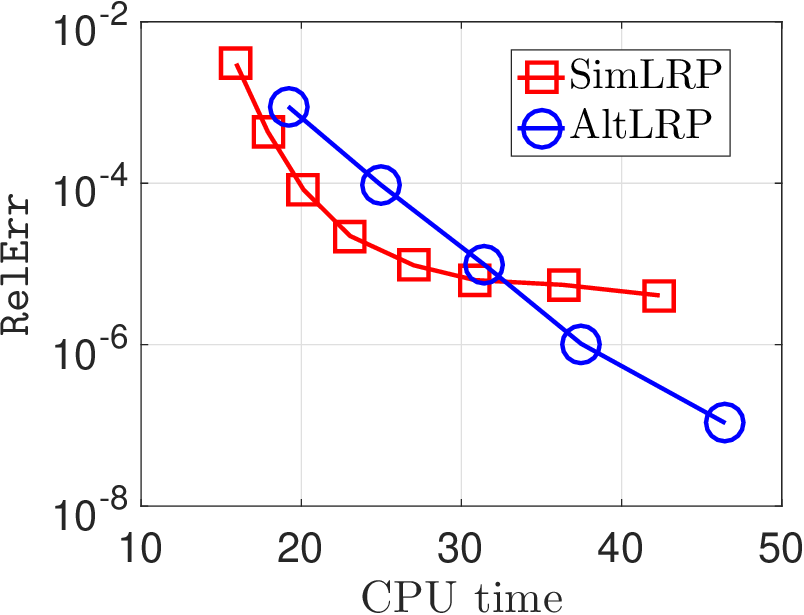}
		}\\
		\subfloat[$N=5, \ \nphy = 257^2$]{
			\includegraphics[width=0.27\linewidth]{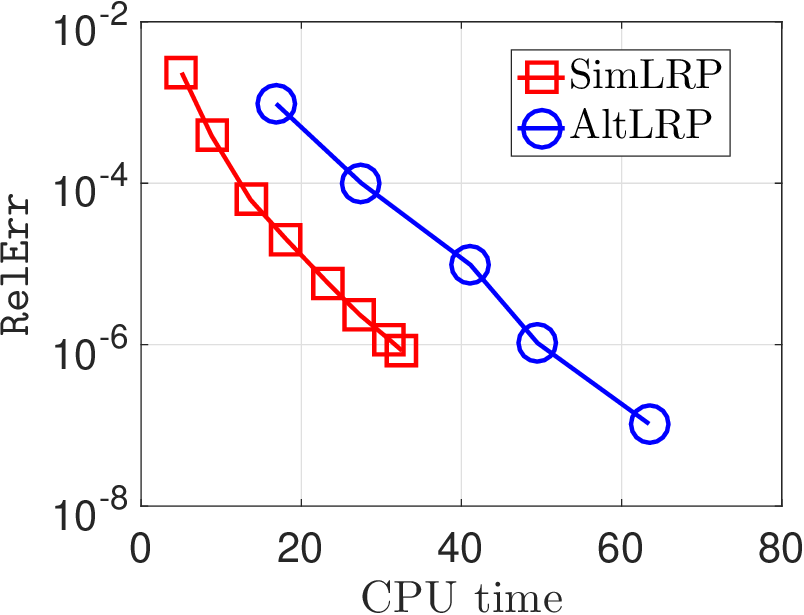}
		}\qquad
		\subfloat[$N=7,\ \nphy = 257^2$]{
			\includegraphics[width=0.27\linewidth]{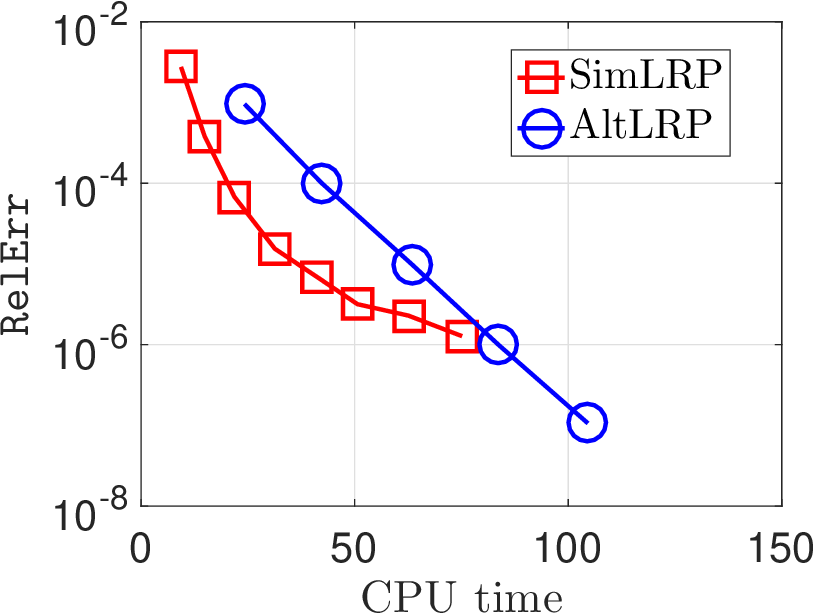}
		}\qquad
		\subfloat[$N=10,\ \nphy = 257^2$]{
			\includegraphics[width=0.27\linewidth]{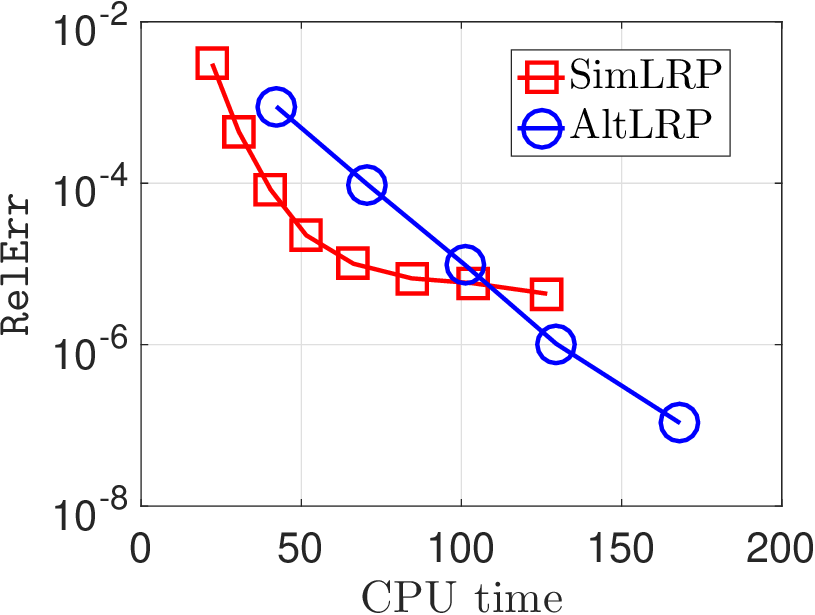}
		}
	\end{center}
	\caption{Relative errors with respect to CPU time in seconds.}\label{fig:test1c}
\end{figure}

In Figure~\ref{fig:test1c}, we present the relative errors with respect to CPU time in seconds. The figure reveals that, for a comparable level of accuracy, the CPU time of the SimLRP approach is slightly smaller than that of the AltLRP approach when $N = 5$ or when the relative error is large. However, for $N=7$ and $N=10$, AltLRP outperforms SimLRP. This observation suggests that, for low-accuracy approximations or low-dimensional problems, the SimLRP approach may be preferred, whereas for higher accuracy or higher dimensional problems, the AltLRP approach is the better choice.

\begin{table}[ht]%
	\caption{CPU time in seconds for standard stochastic Galerkin method.}\label{tab:diff1}
	\begin{center}%
		\newcolumntype{C}{>{\centering\arraybackslash}X}%
		\begin{tabularx}{0.9\linewidth}{CCCCCC}%
			\toprule
			$\nphy$  &   $N=5$ & $N=7$ & $N=10$\\%
			\midrule
			$129^2$     & $38.38$   & $109.37$   &   $394.43$     \\%
			$257^2$     & $229.70$  & $659.54$   &   $2352.73$      \\%
			\bottomrule%
		\end{tabularx}%
	\end{center}%
\end{table}%

Table~\ref{tab:diff1} displays the CPU time in seconds for the standard stochastic Galerkin method. From Figure~\ref{fig:test1c} and Table~\ref{tab:diff1}, it is evident that both the proposed methods are more efficient than the standard stochastic Galerkin method for achieving certain levels of accuracy.

To further evaluate the efficiency of the proposed methods, we compare it with the low-rank solver developed in~\cite{matthies2012solving}, which employs an SVD-based truncation operator $\mathcal{T}_{\tol}$, referred to as LRSSVD in the following discussion. In this test problem, the following truncation strategy is adopted: truncation is applied before preconditioning, after preconditioning, after the computation of the matrix-vector product~(performed in low-rank form), and after every summation during the iterative process of PCG.
    
    \begin{figure}[htbp]
        \begin{center}
            \subfloat[$N=5,\ \nphy = 129^2$]{
                \includegraphics[width=0.27\linewidth]{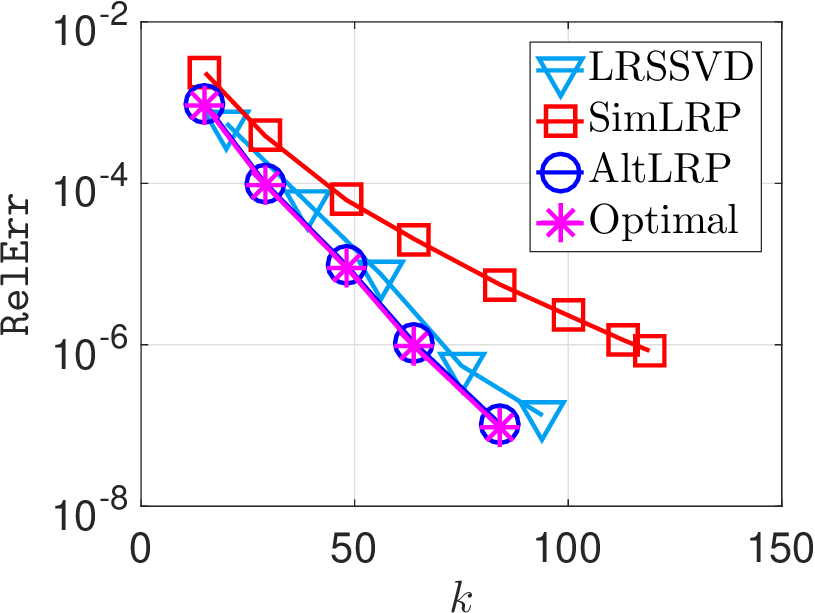}
            }\qquad
            \subfloat[$N=7,\ \nphy = 129^2$]{
                \includegraphics[width=0.27\linewidth]{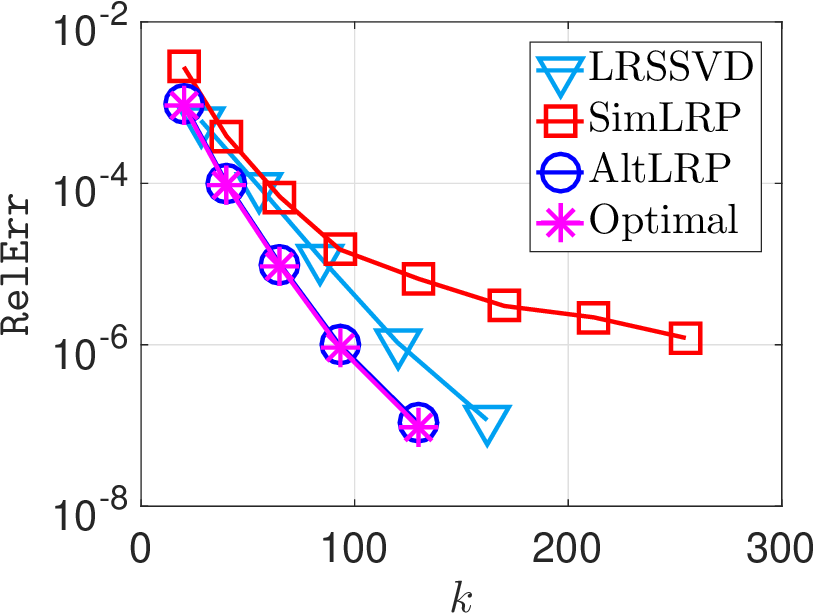}
            }\qquad
            \subfloat[$N=10,\ \nphy = 129^2$]{
                \includegraphics[width=0.27\linewidth]{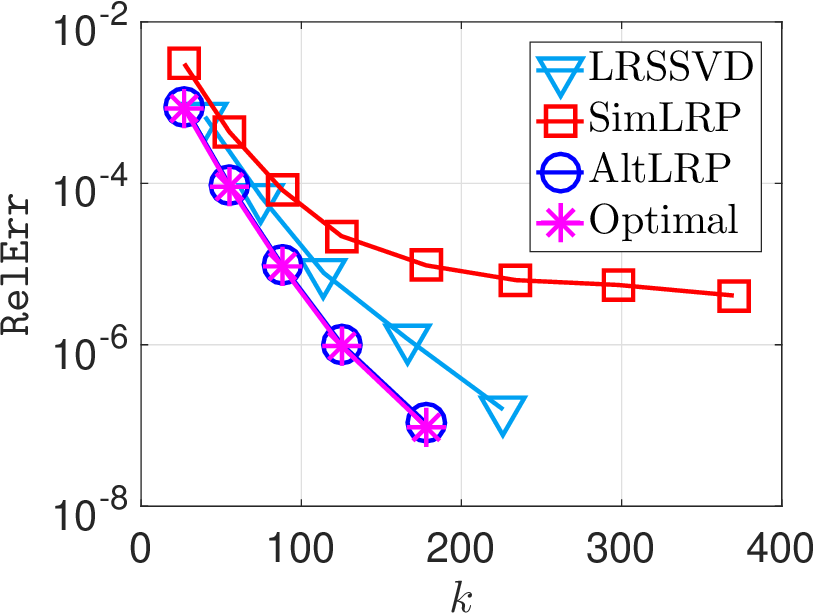}
            }\\
            		\subfloat[$N=5, \ \nphy = 257^2$]{
            			\includegraphics[width=0.27\linewidth]{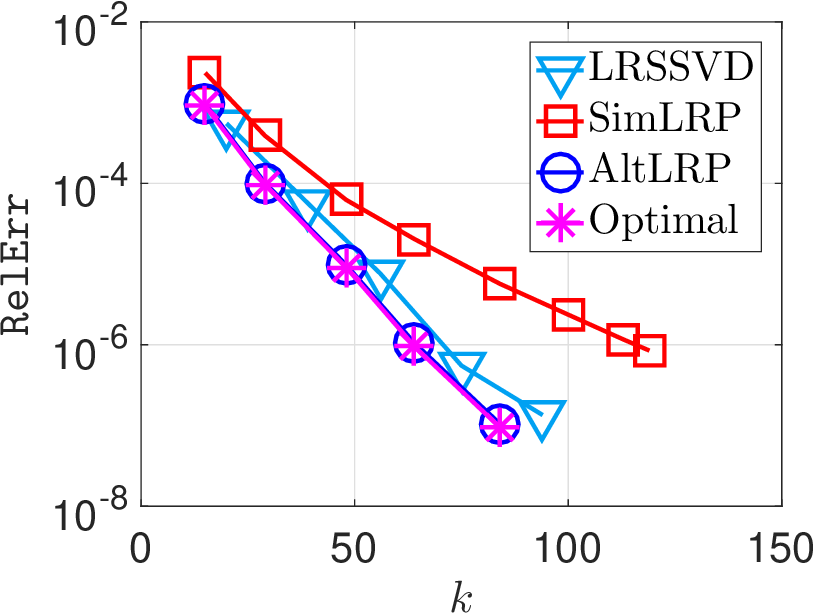}
            		}\qquad
            		\subfloat[$N=7,\ \nphy = 257^2$]{
            			\includegraphics[width=0.27\linewidth]{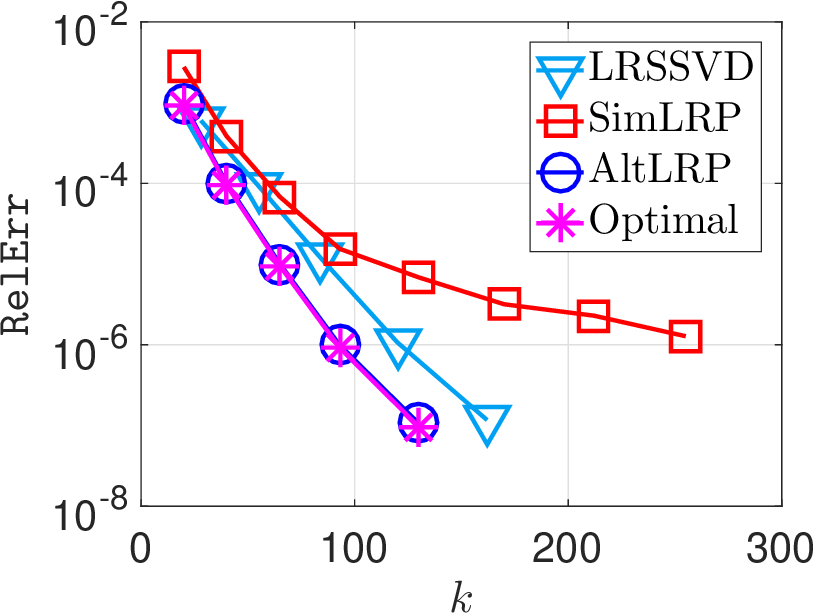}
            		}\qquad
            		\subfloat[$N=10,\ \nphy = 257^2$]{
            			\includegraphics[width=0.27\linewidth]{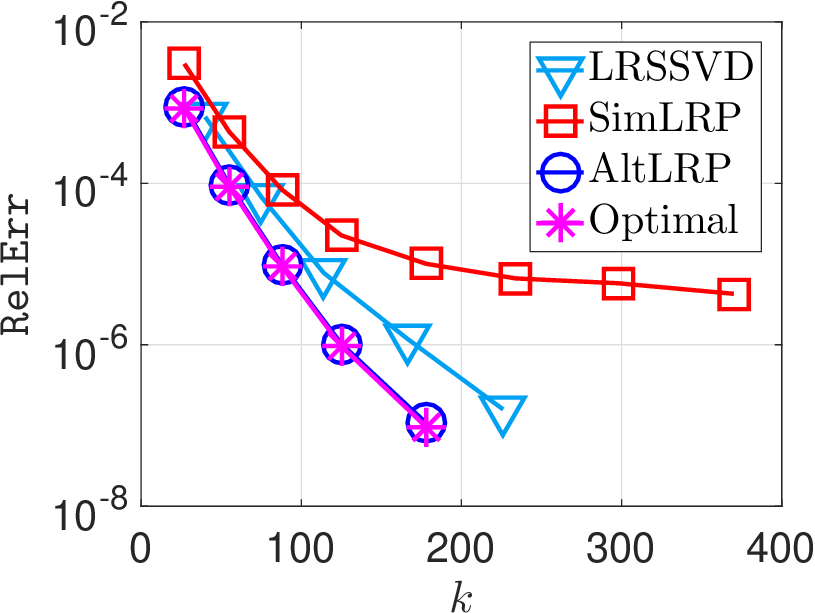}
            		}
        \end{center}
        \caption{Relative errors with respect to rank $k$.}\label{fig:test1d}
    \end{figure}

Figure~\ref{fig:test1d} shows the relative errors with respect to the rank $k$. For the AltLRP approach, the relative error is almost indistinguishable from that of the `optimal' low-rank approximation of the solution. In comparison, the LRSSVD approach requires slightly higher ranks to achieve the same level of accuracy.

The SimLRP approach, however, exhibits a faster growth in rank compared to both the AltLRP and LRSSVD methods---particularly when the stochastic dimension is $N=7$ or $N=10$. In these high-dimensional cases, the SimLRP method shows limited improvement in accuracy as the rank increases, once the relative error reaches approximately $10^{-5}$.

Another noteworthy observation is that, across all methods, the curves of relative error versus rank $k$ remain nearly unchanged for different mesh sizes in the physical space. Note that the relative error of~$\bm{U}$, defined in~\eqref{eq:defRelErr}, can be regarded as a discretized version of the relative error of~$u^{\trunc}(\pv,\sv)$ defined in~\eqref{eq:error_partial sum}. This behavior is, in fact, a direct consequence of Theorem~\ref{thm:quasi_svd}.

\begin{figure}[htbp]
	\begin{center}
		\subfloat[$N=5,\ \nphy = 129^2$]{
			\includegraphics[width=0.27\linewidth]{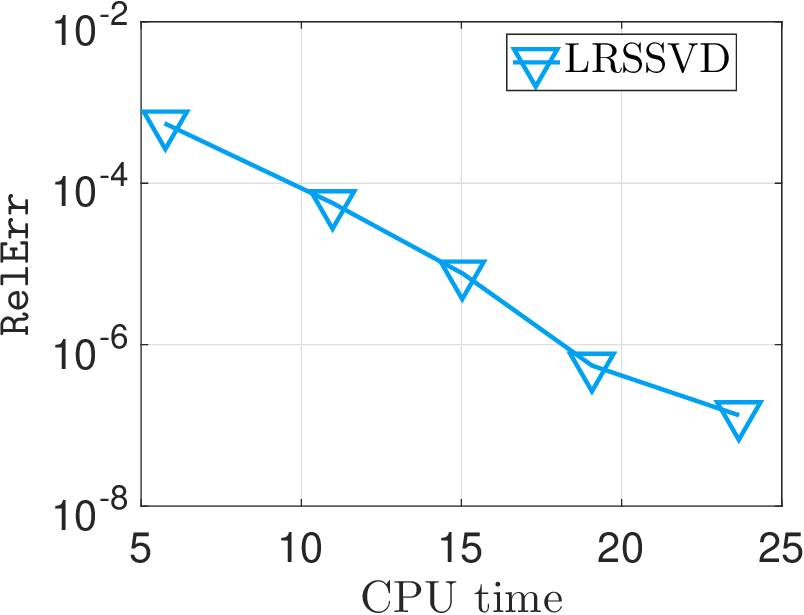}
		}\qquad
		\subfloat[$N=7,\ \nphy = 129^2$]{
			\includegraphics[width=0.27\linewidth]{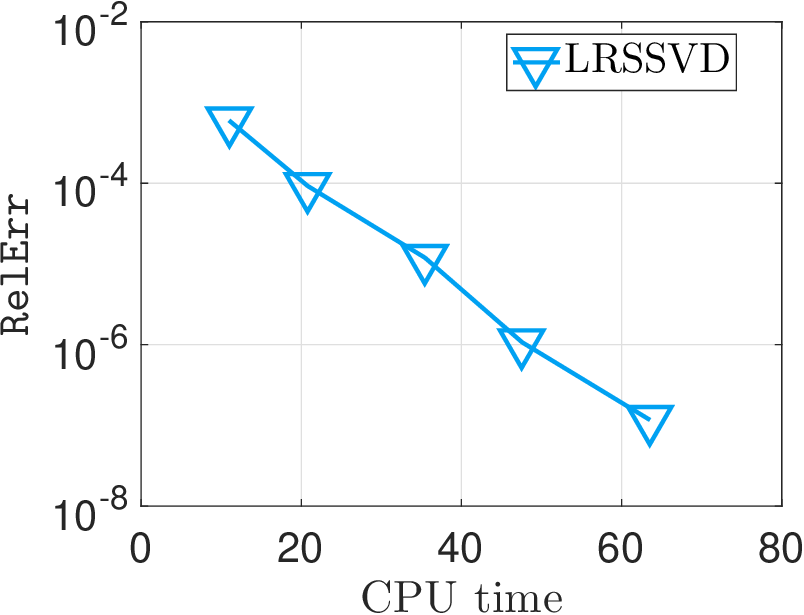}
		}\qquad
		\subfloat[$N=10,\ \nphy = 129^2$]{
			\includegraphics[width=0.27\linewidth]{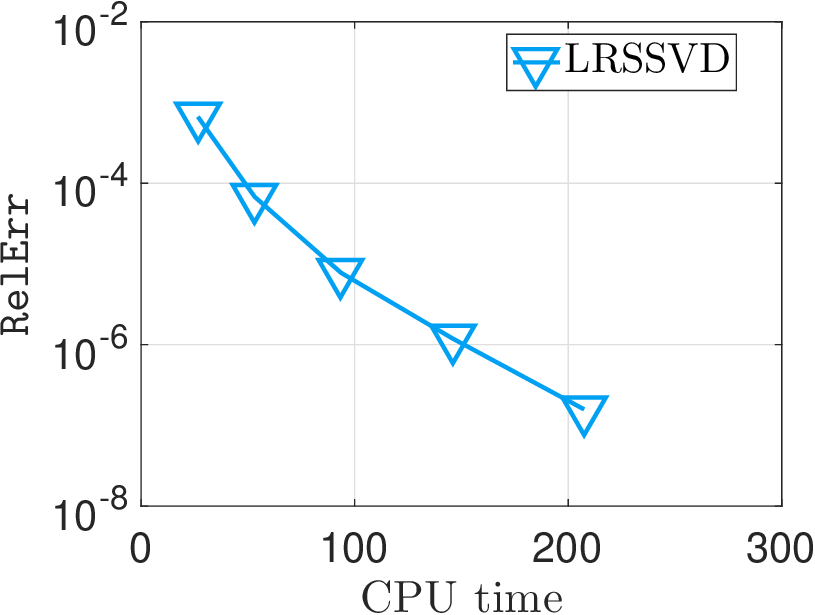}
		}\\
		\subfloat[$N=5, \ \nphy = 257^2$]{
			\includegraphics[width=0.27\linewidth]{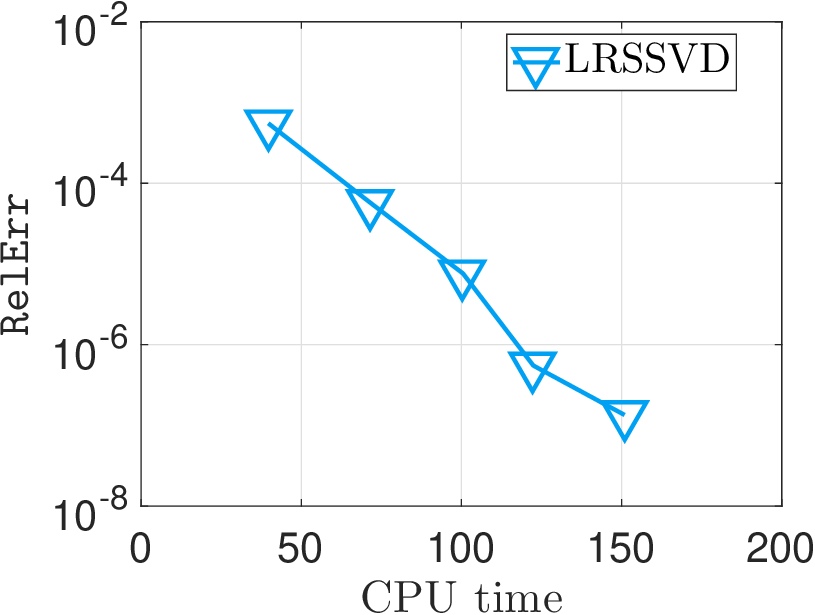}
		}\qquad
		\subfloat[$N=7,\ \nphy = 257^2$]{
			\includegraphics[width=0.27\linewidth]{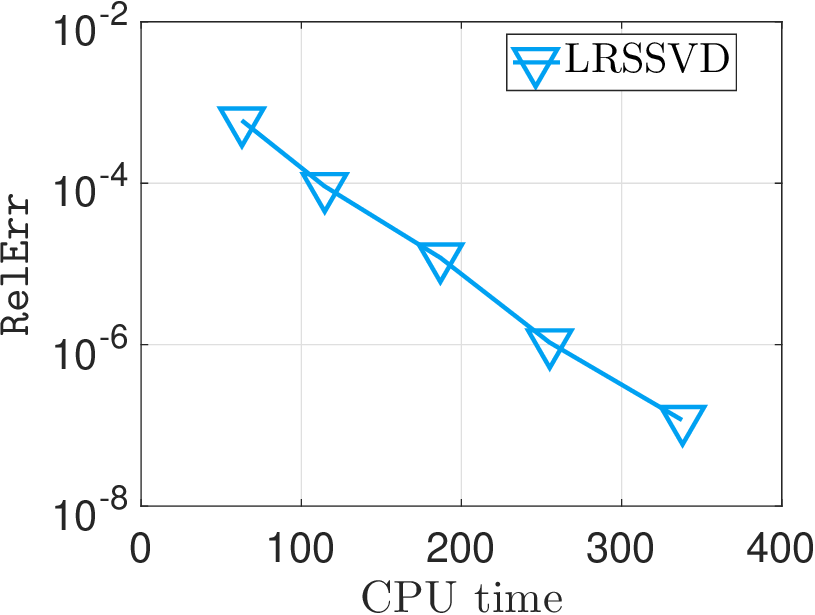}
		}\qquad
		\subfloat[$N=10,\ \nphy = 257^2$]{
			\includegraphics[width=0.27\linewidth]{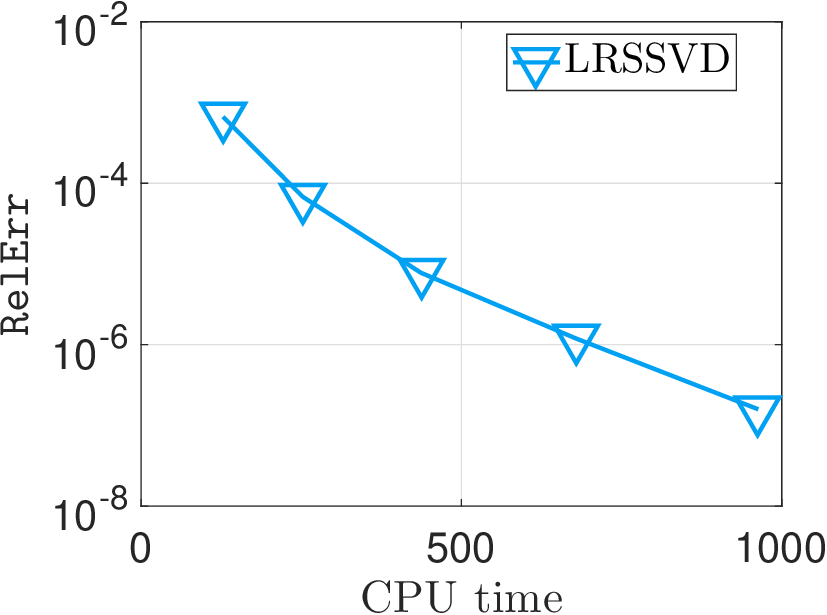}
		}
	\end{center}
	\caption{Relative errors with respect to CPU time in seconds for the LRSSVD approach.}\label{fig:test1e}
\end{figure}

Figure~\ref{fig:test1e} shows the CPU time in seconds for the LRSSVD approach. From Figure~\ref{fig:test1c} and~\ref{fig:test1e}, it is clear that the AltLRP approach requires significantly less CPU time than the LRSSVD approach to achieve a given accuracy. In addition, if the required relative error is above $10^{-5}$, the SimLRP approach can also be applied.

The efficiency of the AltLRP approach stems from two main reasons. First, the method proposed in this work determines the rank required for a given accuracy 
by analyzing the singular values of quasimatrices associated with the truncated gPC expansion of the solution on a coarse physical grid, allowing the construction of a fixed-rank approximation. This avoids the need for truncation during the PCG iterations for solving the linear systems. Second, while the LRSSVD method operates directly on the full linear system and applies low-rank techniques during the iteration, the proposed AltLRP approach constructs a much smaller projected linear system based on the identified low-rank subspace. This projection significantly reduces computational cost compared to the LRSSVD approach.

It should be pointed out that the computational time in LRSSVD depends on several factors, such as the method used for SVD updates and the termination criteria of the iteration. In this work, we only employ a simple implementation of ~\cite{matthies2012solving}, which may not achieve its optimal efficiency.

\subsection{Test problem 2}

In this test problem, we consider the stochastic Helmholtz problem given by
\begin{equation}\notag
\nabla^2u + a^2(\pv,\sv)u=f(\pv) \quad \mbox{in}\quad D\times\Gamma,
\end{equation}
with the Sommerfeld radiation boundary condition.
Here, $D = [0,1]^2$ is the domain of interest, and the Helmholtz coefficient $a(\pv,\sv)$ is a truncated KL expansion of a random field with a mean function  $a_0(\pv)=4\cdot(2\pi)$, a standard deviation $\sigma=0.8\pi$, and the covariance function $\mathrm{Cov}\,(\bm{x},\bm{y})$ given by
\begin{equation}\notag
\mathrm{Cov}\,(\bm{x},\bm{y})=\sigma^2 \exp\left(-\frac{|x_1-y_1|}{c}-\frac{|x_2-y_2|}{c}
\right),
\end{equation}
where $\pv=[x_1,x_2]^T$, $\bm{y}=[y_1,y_2]^T$ and $c=4$ is the correlation length. Note that the KL expansion takes the form
\begin{equation}\label{eq:kl2}\notag
a(\pv,\sv)=a_0(\pv)+\sum_{i=1}^{N}a_i(\pv)\svc_i=a_0(\pv)+\sum_{i=1}^{N}\sqrt{\lambda_i}c_i(\pv)\svc_i,
\end{equation}
where $\{\lambda_i,c_i(\pv)\}_{i=1}^{N}$ are the eigenpairs of  $\mathrm{Cov}\,(\pv,\bm{y})$,
$\{\svc_i\}^{N}_{i=1}$ are uncorrelated random variables, and $N$ is the number of KL modes retained. The  Gaussian point source at the center of the domain is used as the source term, i.e.,
\begin{equation}\notag
f(\pv) = \mbox{e}^{-(8\cdot4)^2((x_1-0.5)^2+(x_2-0.5)^2)}.
\end{equation}

For this test problem, we assume that the random variables $\{\svc_i\}^{N}_{i=1}$ are independent and uniformly distributed within the range $[-1,1]$. We use the perfectly matched layers (PML) to simulate the Sommerfeld condition~\cite{Berenger94}, and discretize the physical space using $Q_1$ finite elements based on the IFISS package~\cite{ifiss-siamreview}.

In this test problem, the gPC order $p$ is set to $4$, the degrees of freedom (DOF) of the coarse grid in physical space (i.e., $\bm{T}^c$) is set to $\nphy = 33\times 33$, and the cardinality of the random set~$\bm{\Theta}$ is set to $|\bm{\Theta}| = \min(2\nstoch,\nphy)$. All the linear systems,~i.e.,~\eqref{eq:SimLRP} and~\eqref{eq:AltLRP_stoch}--\eqref{eq:AltLRP_phy}, are solved using the preconditioned bi-conjugate gradient stabilized~(Bi-CGSTAB) method with a mean-based preconditioner~\cite{Powell2009} and a tolerance of $10^{-8}$.

\begin{figure}[htbp]
	\begin{center}
		\subfloat[$N=5$]{
			\includegraphics[width=0.27\linewidth]{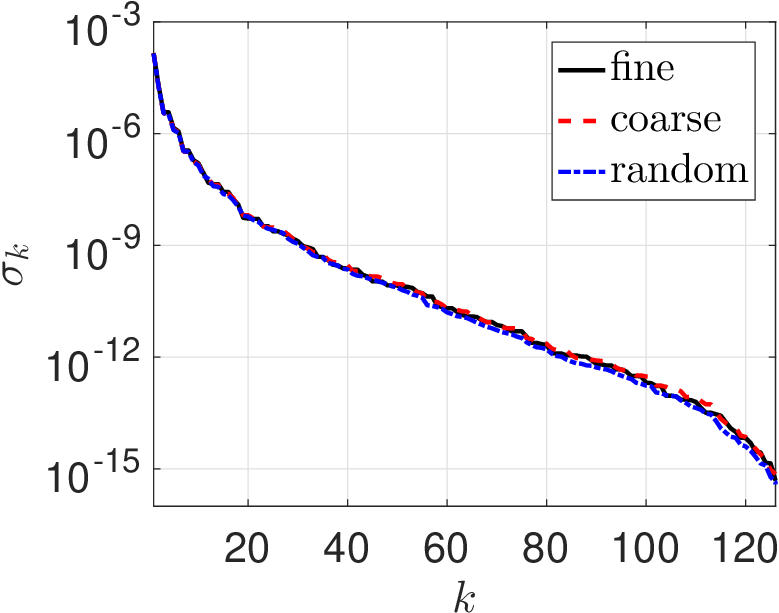}
		}\qquad
		\subfloat[$N=7$]{
			\includegraphics[width=0.27\linewidth]{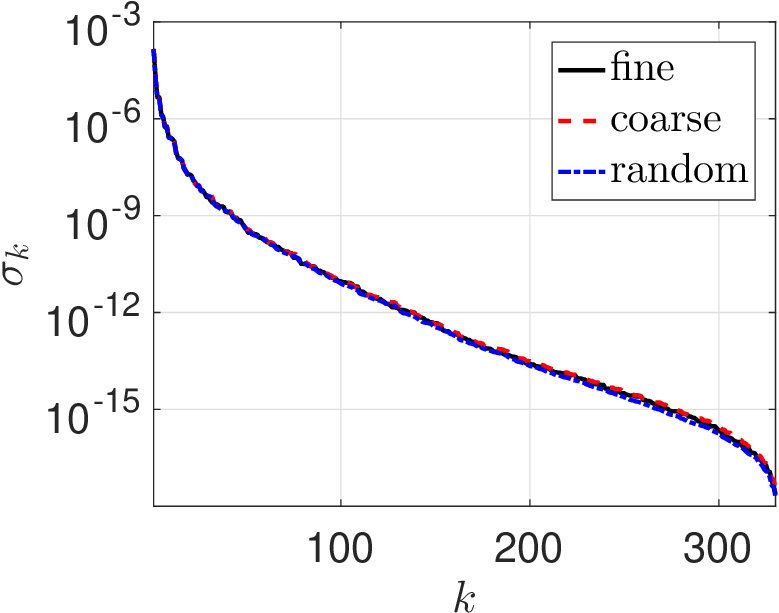}
		}\qquad
		\subfloat[$N=10$]{
			\includegraphics[width=0.27\linewidth]{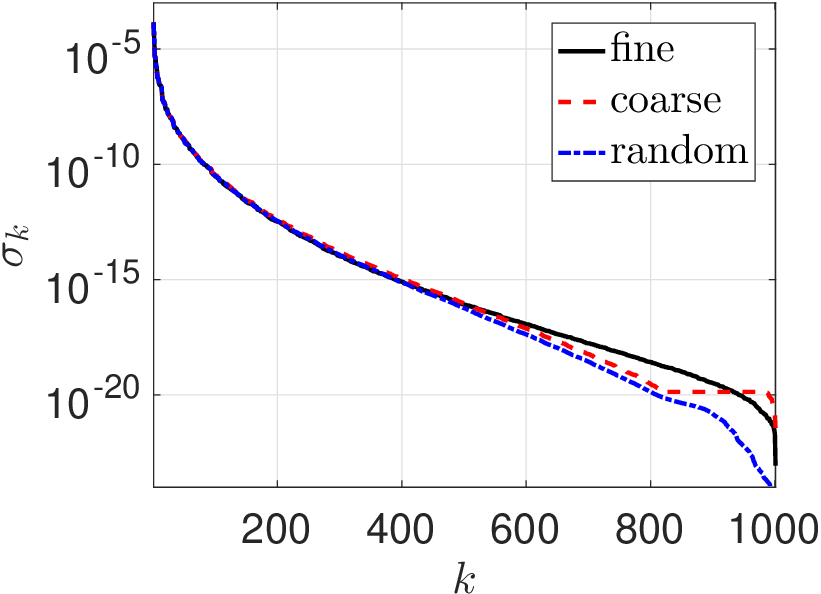}
		}
	\end{center}
	
	\caption{Singular values for different $m$ and $\tol$.}\label{fig:test2SVD}
	{\small In this figure, `fine' refers to the singular values of $\quasimatrix{U}^{\trunc}_{D\times \nstoch}$ with $\nphy= 257^2$; `coarse' refers to the singular values of $\quasimatrix{U}^{\coarse}_{D\times \nstoch}$ with $\nphy= 33^2$;  and `random' refers to the singular values of  $\widehat{\quasimatrix{U}}^{\samp}_{D\times n}$ in Corollary~\ref{co:snap_equv} with $|\bm{\Theta}|=\min(2\nstoch,\nphy)$ and $\nphy=33^2$.}
\end{figure}

In Figure~\ref{fig:test2SVD}, we plot the singular values for different $N$. The DOF of the physical space is $\nphy = 33^2$ for the coarse grid and $\nphy = 257^2$ for the fine grid. It is noteworthy that, in this example, the $L^2$ norm of a function is simply approximated by  the $2$-norm of its corresponding vector, multiplied by a factor of $1/\nphy$, since the mesh grid in physical space is uniform.  From the figure, we observe that the singular values for `fine', `coarse' and `random' coincide very well with each other, reflecting the conclusion of Theorem~\ref{th:snap_equv}. It also demonstrates that the singular values of~$\quasimatrix{U}^{\trunc}_{D\times \nstoch}$ are not very sensitive to the DOF of the physical mesh grid, and thus allowing the use of singular values computed on a coarse grid to determine the rank $k$ for a desired accuracy $\tol$.

\begin{figure}[htbp]	
	\begin{center}
		\subfloat[$N=5,\ \nphy = 129^2$]{
			\includegraphics[width=0.27\linewidth]{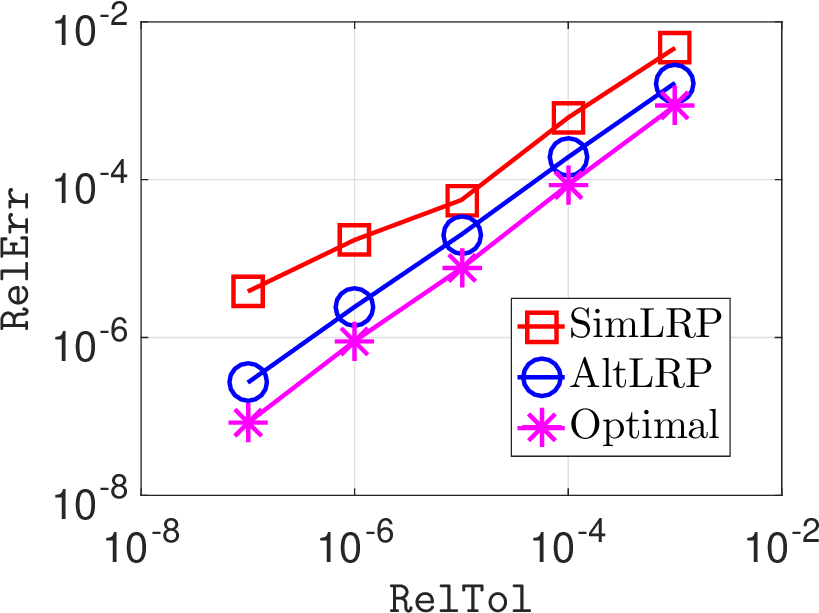}
		}\qquad
		\subfloat[$N=7,\ \nphy = 129^2$]{
			\includegraphics[width=0.27\linewidth]{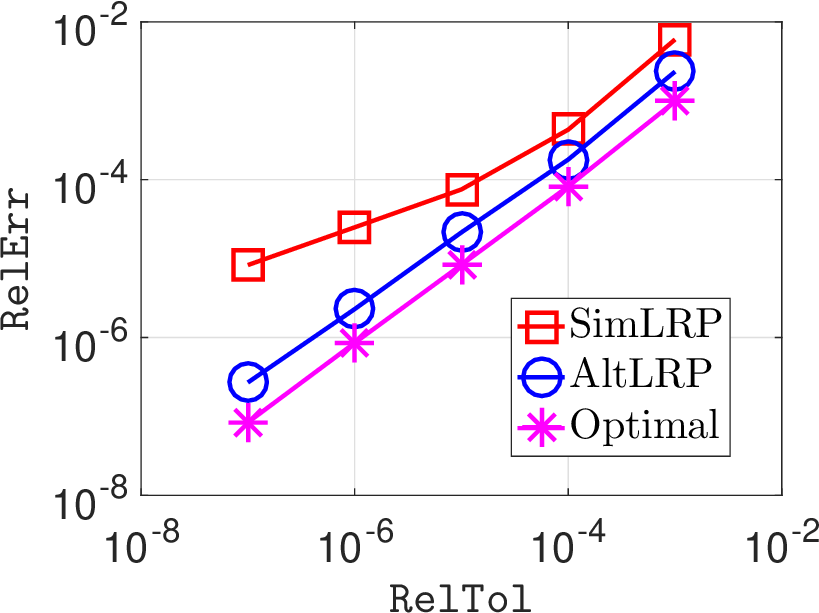}
		}\qquad
		\subfloat[$N=10,\ \nphy = 129^2$]{
			\includegraphics[width=0.27\linewidth]{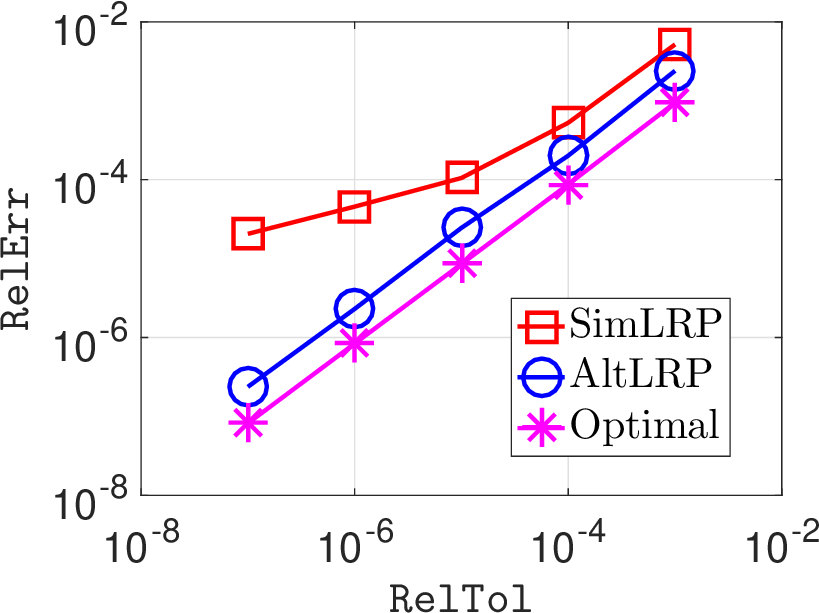}
		}\\
		\subfloat[$N=5, \ \nphy = 257^2$]{
			\includegraphics[width=0.27\linewidth]{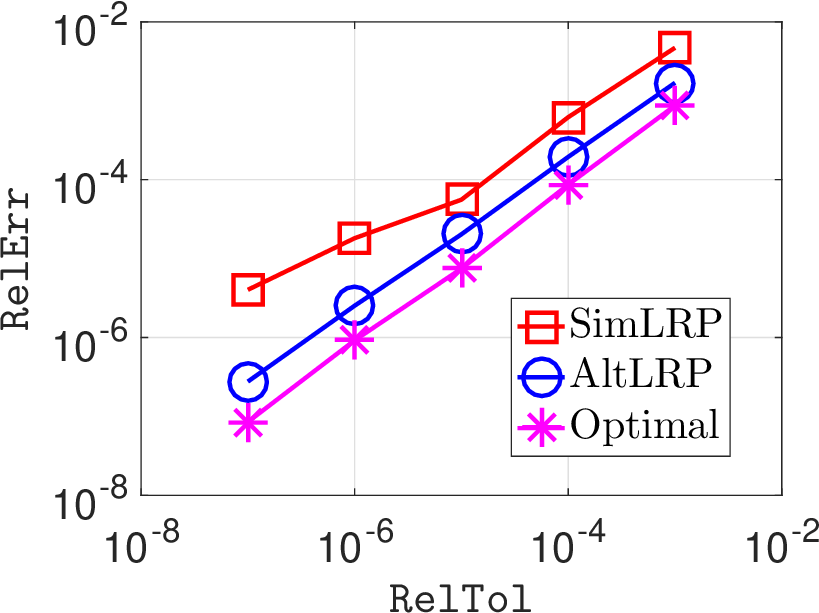}
		}\qquad
		\subfloat[$N=7,\ \nphy = 257^2$]{
			\includegraphics[width=0.27\linewidth]{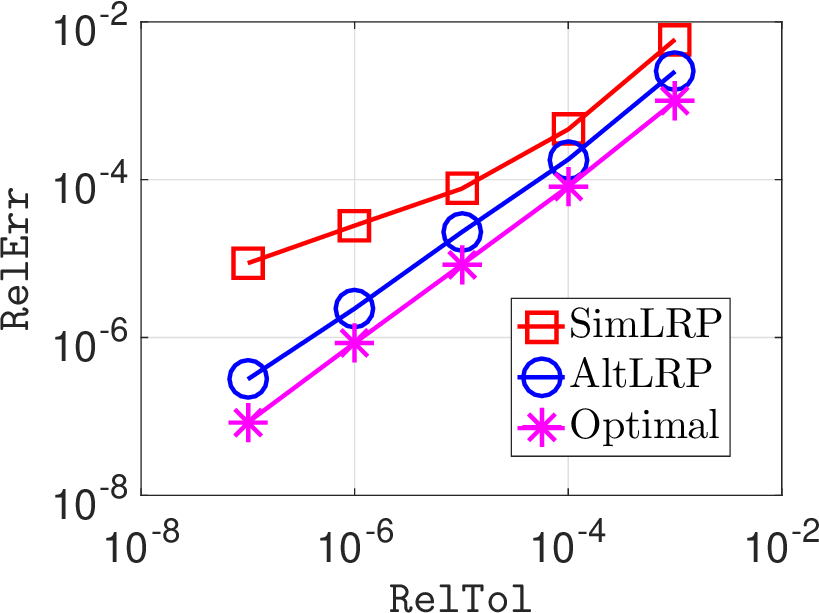}
		}\qquad
		\subfloat[$N=10,\ \nphy = 257^2$]{
			\includegraphics[width=0.27\linewidth]{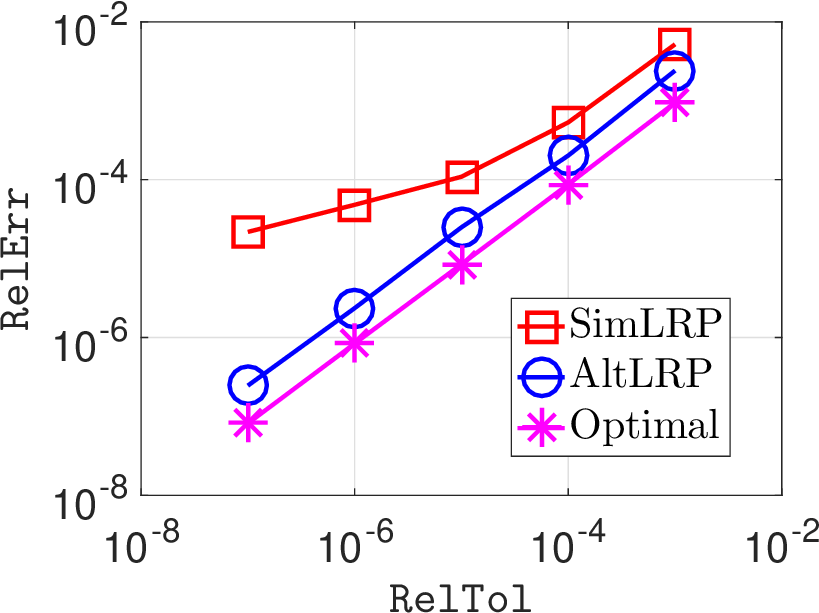}
		}
	\end{center}
	\caption{Relative errors with respect to different $\tol$.}\label{fig:test2b}
	
\end{figure}

In Figure~\ref{fig:test2b}, we present the relative errors concerning different values of $\tol$, where the relative error is defined by~\eqref{eq:defRelErr}. The `optimal' low-rank approximation of the solution is obtained through the SVD decomposition of $\bm{U}$ with the corresponding $k$ terms retained. In the AltLRP approach, the value of $i_{\max}$ is set to be~$1$. From the figure, it is evident that the AltLRP approach performs nearly optimally. While the SimLRP approach also demonstrates good performance, the tolerance indicator~\eqref{eq:find_k} for it deteriorates for smaller values of $\tol$.

\begin{figure}[htbp]
	\begin{center}
		\subfloat[$N=5,\ \nphy = 129^2$]{
			\includegraphics[width=0.27\linewidth]{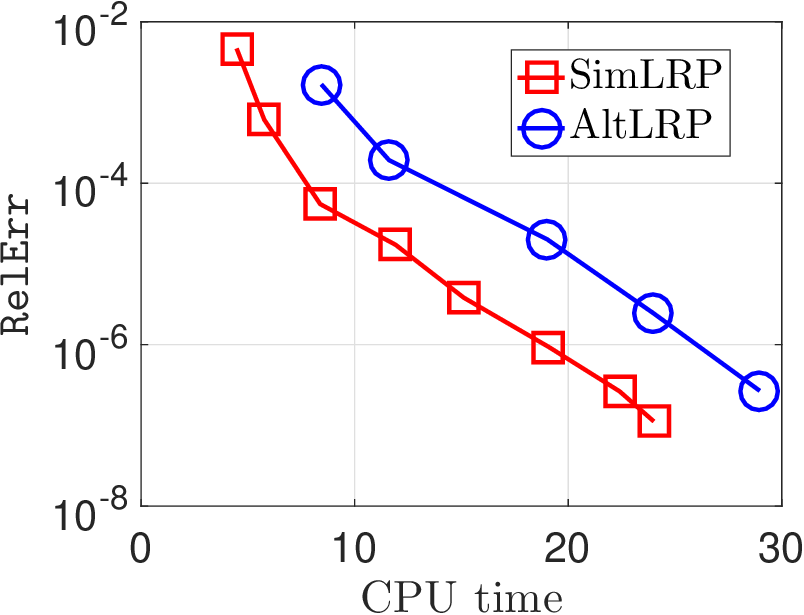}
		}\qquad
		\subfloat[$N=7,\ \nphy = 129^2$]{
			\includegraphics[width=0.27\linewidth]{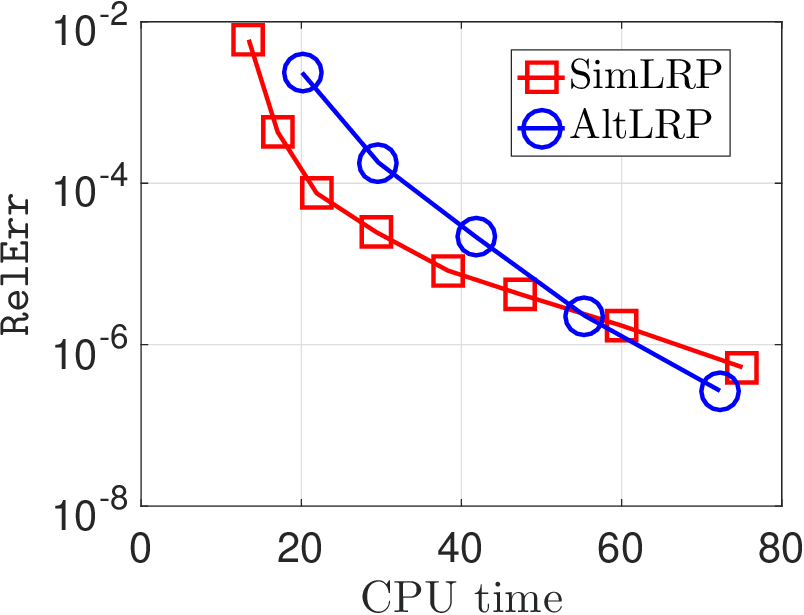}
		}\qquad
		\subfloat[$N=10,\ \nphy = 129^2$]{
			\includegraphics[width=0.27\linewidth]{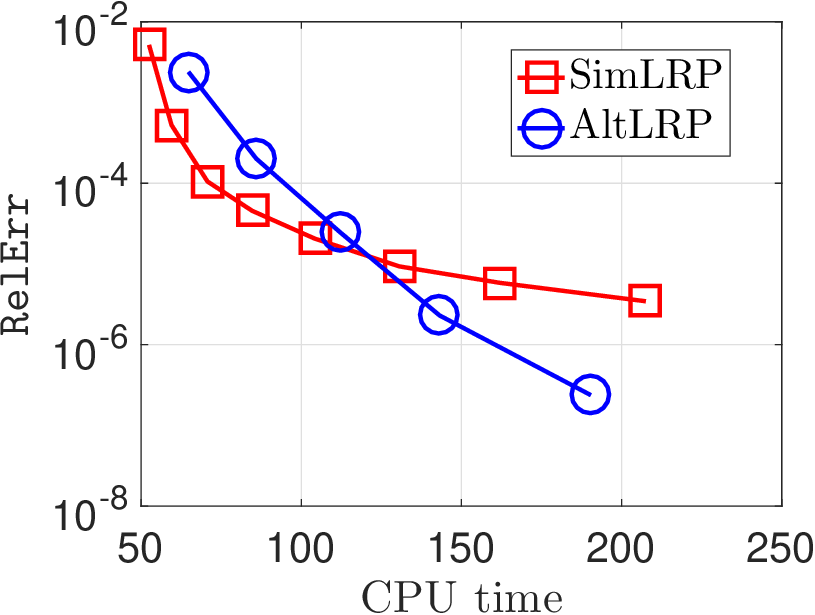}
		}\\
		\subfloat[$N=5, \ \nphy = 257^2$]{
			\includegraphics[width=0.27\linewidth]{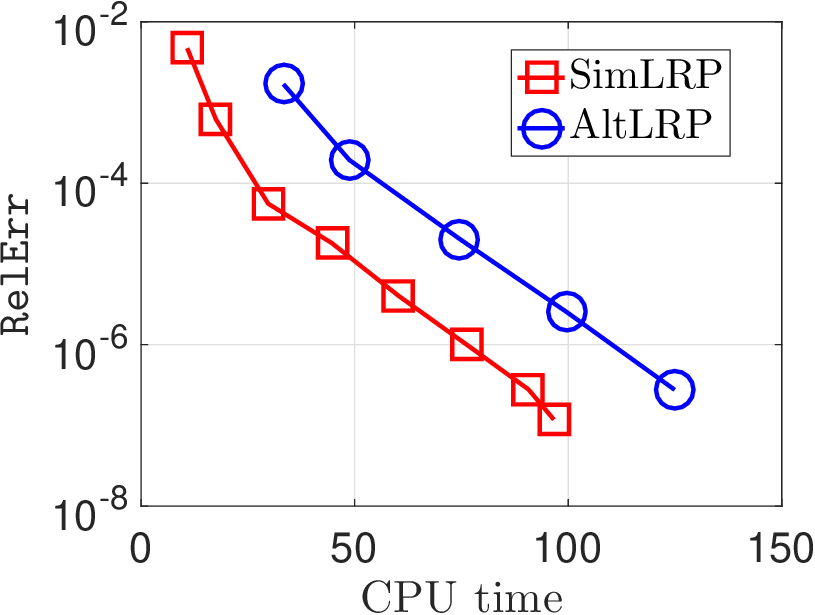}
		}\qquad
		\subfloat[$N=7,\ \nphy = 257^2$]{
			\includegraphics[width=0.27\linewidth]{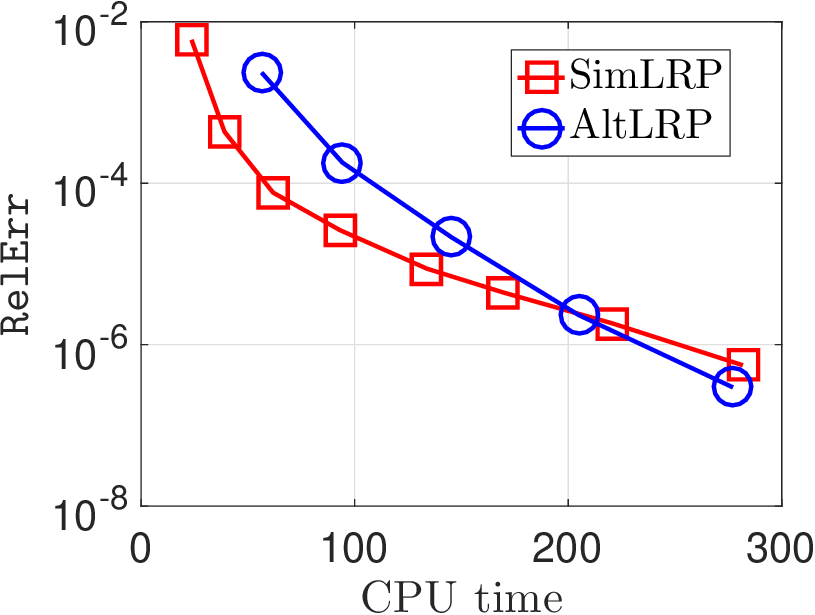}
		}\qquad
		\subfloat[$N=10,\ \nphy = 257^2$]{
			\includegraphics[width=0.27\linewidth]{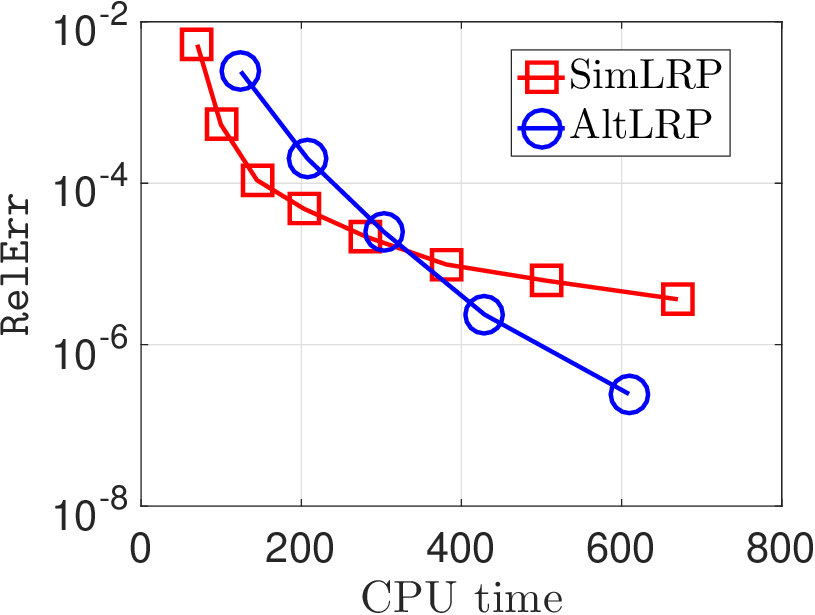}
		}
	\end{center}
	\caption{Relative errors with respect to CPU time in seconds.}\label{fig:test2c}
	
\end{figure}

In Figure~\ref{fig:test2c}, we present the relative errors with respect to CPU time in seconds. The figure reveals that, for a comparable level of accuracy, the CPU time of the SimLRP approach is slightly smaller than that of the AltLRP approach when $N = 5$ or when the relative error is large. However, for $N=7$ and $N=10$, AltLRP outperforms SimLRP. This observation suggests that, for low-accuracy approximations or low-dimensional problems, the SimLRP approach may be preferred, whereas for higher accuracy or higher dimensional problems, the AltLRP approach is the better choice.

\begin{table}[ht]%
	\caption{CPU time in seconds for standard stochastic Galerkin method.}\label{tab:helm1}
	\begin{center}%
		\newcolumntype{C}{>{\centering\arraybackslash}X}%
		\begin{tabularx}{0.9\linewidth}{CCCCCC}%
			\toprule
			$\nphy$  	&   $N=5$ & $N=7$ & $N=10$\\%
			\midrule
			$129^2$     &  $57.74$  & $214.71$    &  $908.00$      \\%
			$257^2$     &  $308.87$ & $1138.43$   &  $4645.46$      \\%
			\bottomrule%
		\end{tabularx}%
	\end{center}%
\end{table}%

Table~\ref{tab:helm1} displays the CPU time in seconds for the standard stochastic Galerkin method. Again, from Figure~\ref{fig:test2c} and Table~\ref{tab:helm1}, it is evident that both the proposed methods are more efficient than the standard stochastic Galerkin method for achieving certain levels of accuracy.

To further evaluate the efficiency of the proposed method, we compare it with the LRSSVD approach developed in~\cite{matthies2012solving}. The truncation strategy used here is the same as in the previous test problem: truncation is applied before preconditioning, after preconditioning, after the matrix-vector product (performed in low-rank form), and after each summation during the iterative process of Bi-CGSTAB.

It is worth noting that although this test problem involves only $N$  random variables, the number of terms on the left-hand side of equation~\eqref{eq:linsg_matrix} is $K=(N+1)^2+1$, which is significantly larger than that in Test Problem 1. This leads to increased memory requirements during the iterative process, since the rank after the matrix-vector product (in low-rank form) becomes $K\times k$. As a result, some cases run out of memory.

\begin{figure}[htbp]
	\begin{center}
		\subfloat[$N=5,\ \nphy = 129^2$]{
			\includegraphics[width=0.27\linewidth]{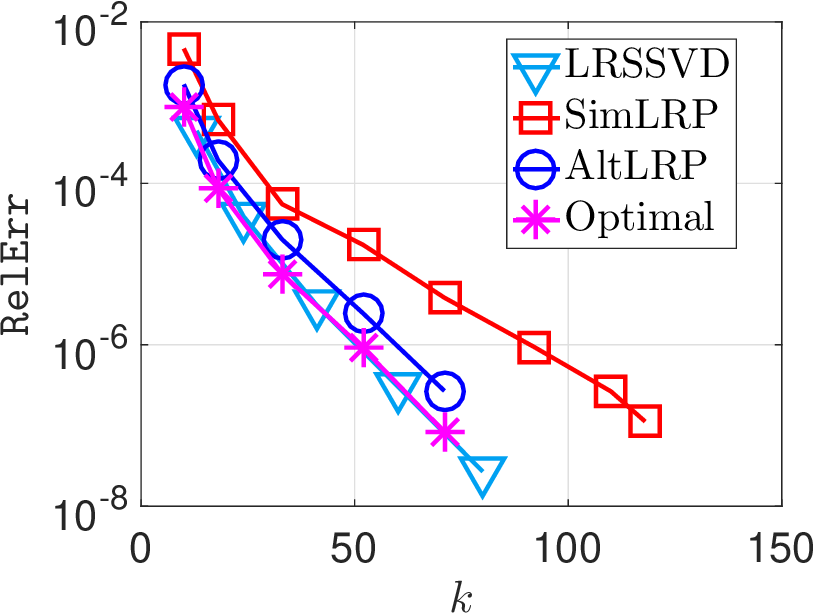}
		}\qquad
		\subfloat[$N=7,\ \nphy = 129^2$]{
			\includegraphics[width=0.27\linewidth]{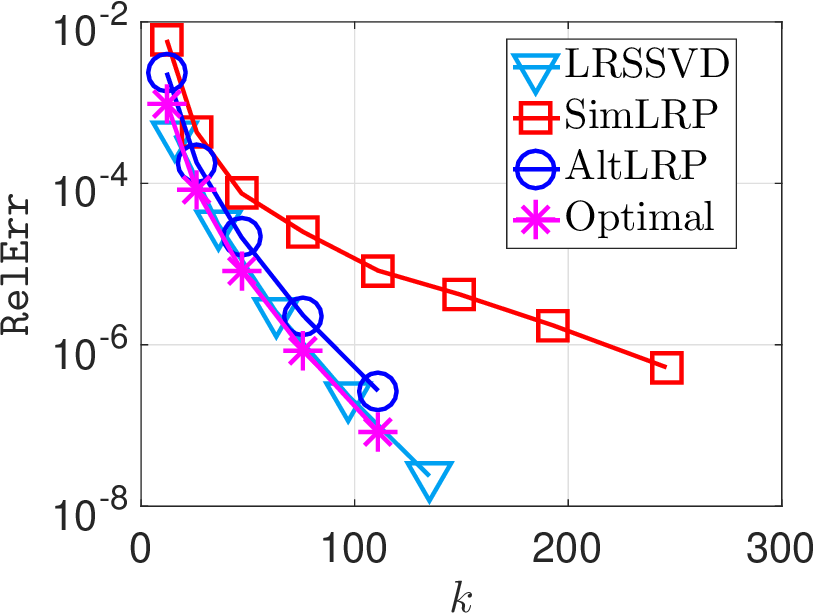}
		}\qquad
		\subfloat[$N=10,\ \nphy = 129^2$]{
			\includegraphics[width=0.27\linewidth]{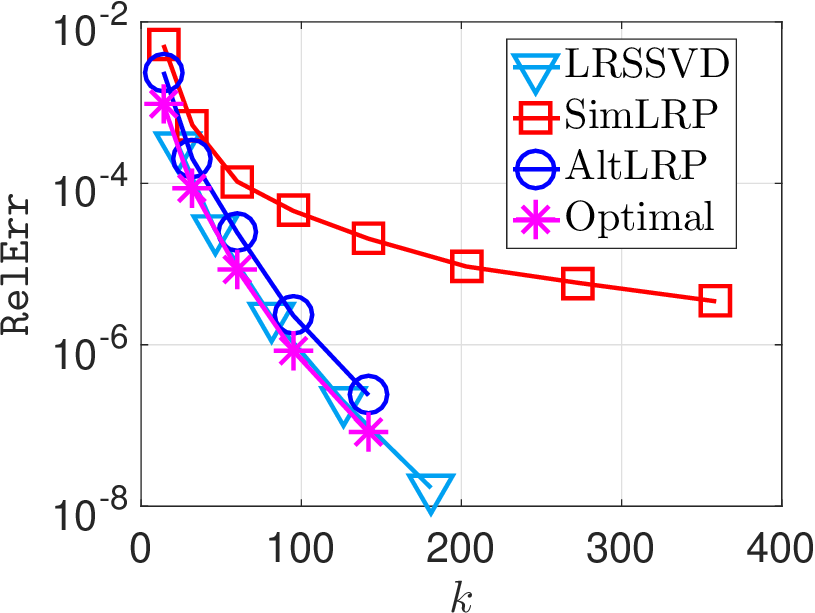}
		}\\
		\subfloat[$N=5, \ \nphy = 257^2$]{
			\includegraphics[width=0.27\linewidth]{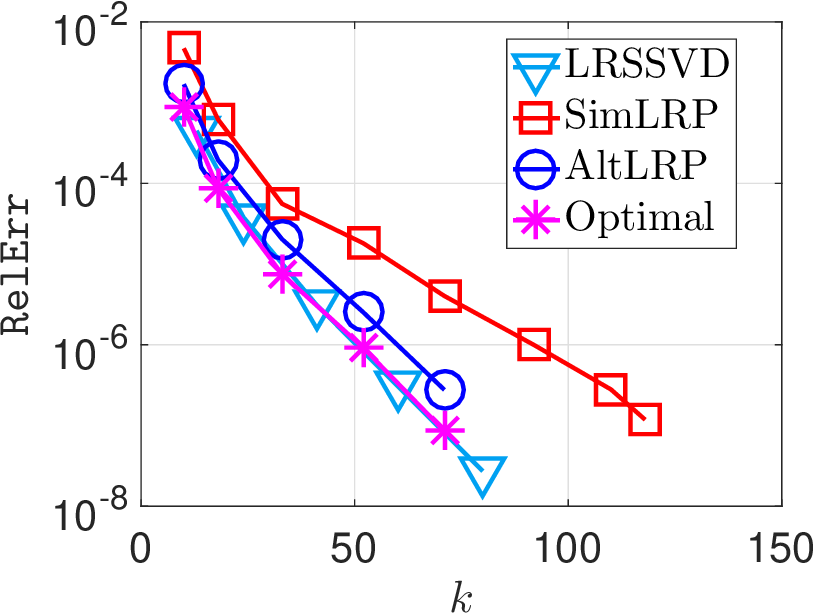}
		}\qquad
		\subfloat[$N=7,\ \nphy = 257^2$]{
			\includegraphics[width=0.27\linewidth]{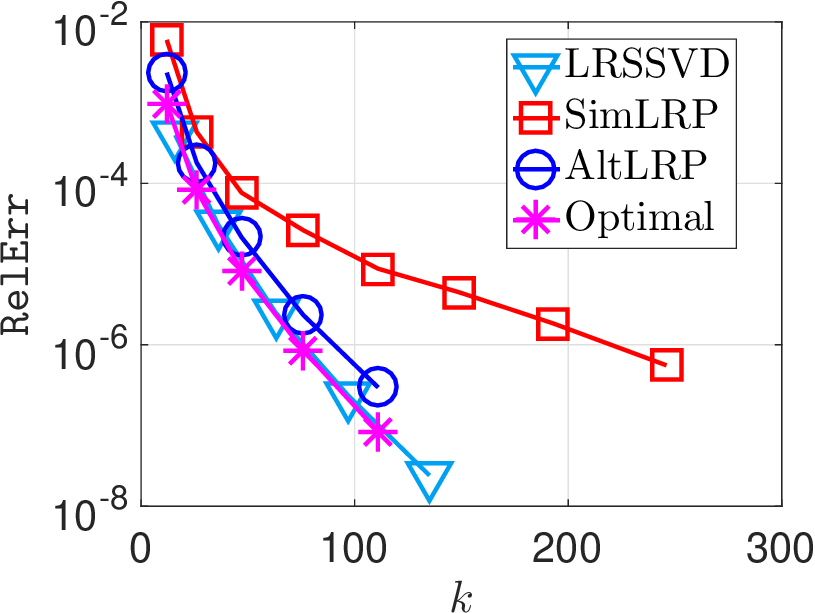}
		}\qquad
		\subfloat[$N=10,\ \nphy = 257^2$]{
			\includegraphics[width=0.27\linewidth]{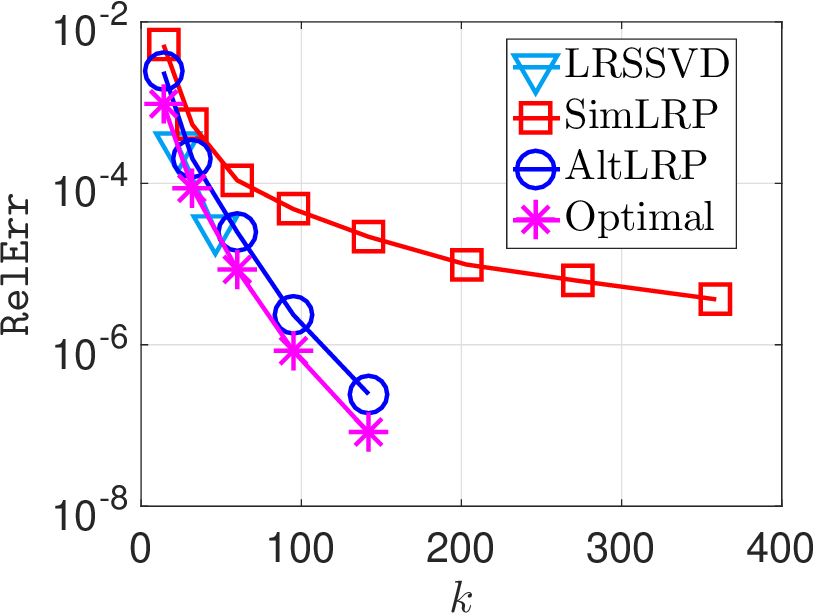}
		}
	\end{center}
	\caption{Relative errors with respect to rank $k$.}\label{fig:test2d}
\end{figure}

Figure~\ref{fig:test2d} shows the relative errors with respect to the rank $k$. For the LRSSVD approach, the relative error is almost indistinguishable from that of the `optimal' low-rank approximation of the solution. In comparison, the AltLRP approach requires slightly higher ranks to achieve the same level of accuracy. The SimLRP approach, however, exhibits a faster growth in rank compared to both the AltLRP and LRSSVD methods---particularly when the stochastic dimension is $N=7$ or $N=10$. In these high-dimensional cases, the SimLRP method shows limited improvement in accuracy as the rank increases, once the relative error reaches approximately $10^{-5}$.

Another noteworthy observation is that, across all methods, the curves of relative error versus rank $k$ remain nearly unchanged for different mesh sizes in the physical space. Note that the relative error of~$\bm{U}$, defined in~\eqref{eq:defRelErr}, can be regarded as a discretized version of the relative error of~$u^{\trunc}(\pv,\sv)$ defined in~\eqref{eq:error_partial sum}. This behavior is, in fact, a direct consequence of Theorem~\ref{thm:quasi_svd}.

\begin{figure}[htbp]
	\begin{center}
		\subfloat[$N=5,\ \nphy = 129^2$]{
			\includegraphics[width=0.27\linewidth]{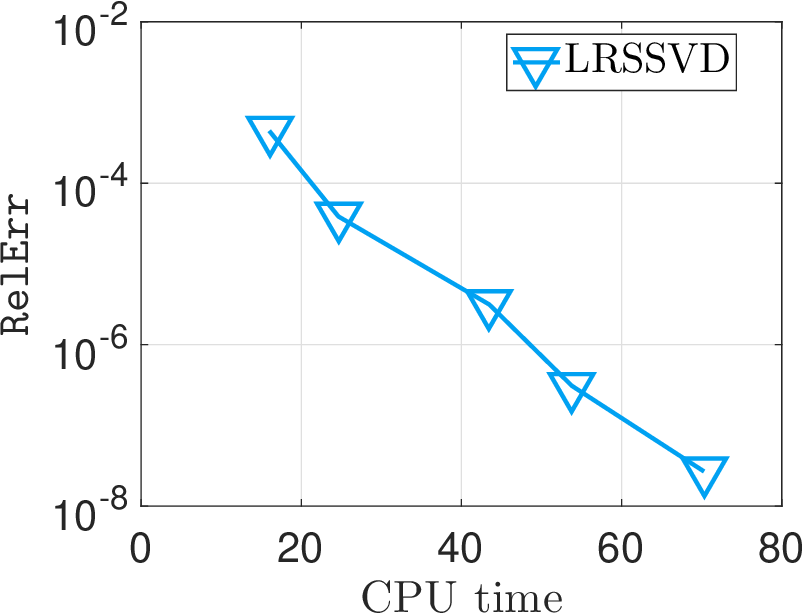}
		}\qquad
		\subfloat[$N=7,\ \nphy = 129^2$]{
			\includegraphics[width=0.27\linewidth]{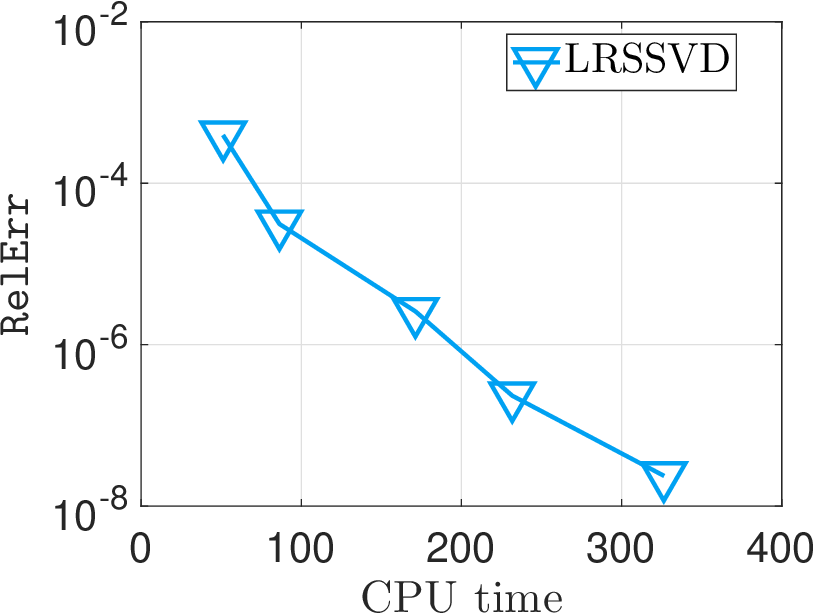}
		}\qquad
		\subfloat[$N=10,\ \nphy = 129^2$]{
			\includegraphics[width=0.27\linewidth]{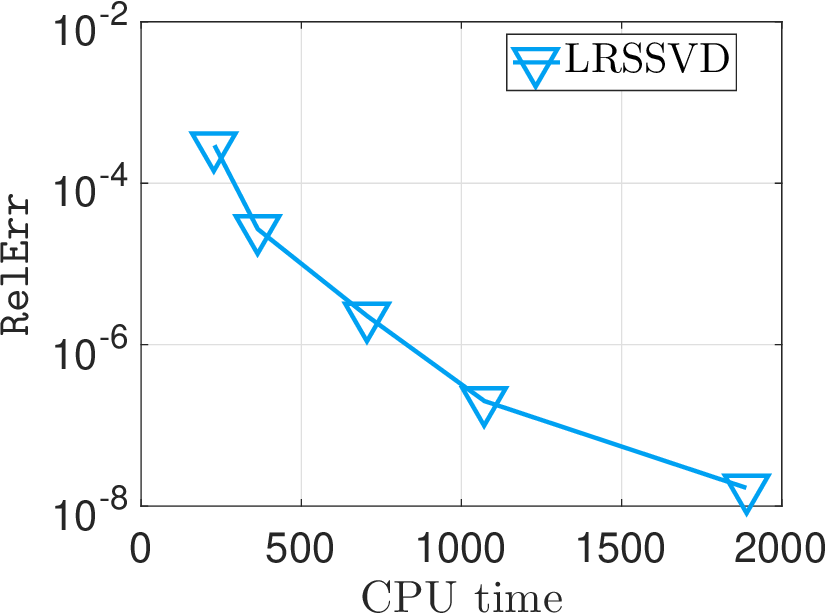}
		}\\
		\subfloat[$N=5, \ \nphy = 257^2$]{
			\includegraphics[width=0.27\linewidth]{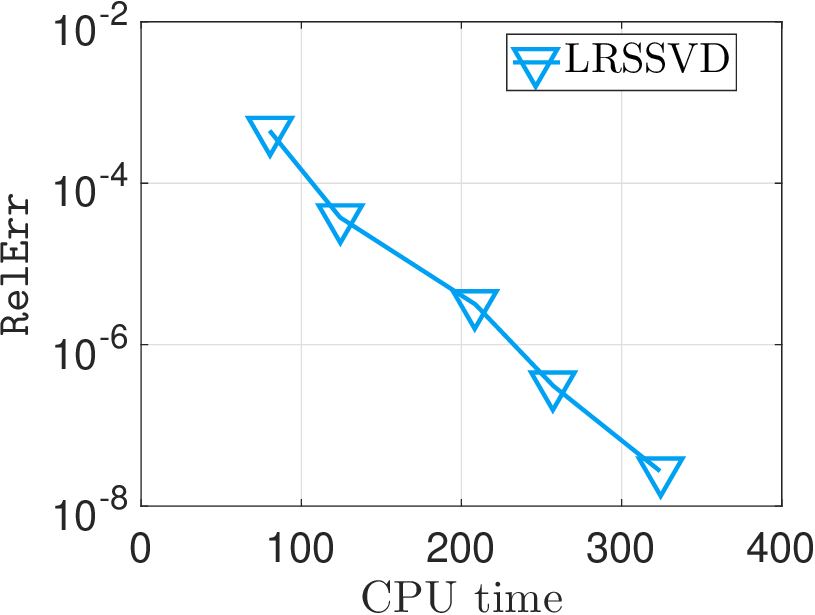}
		}\qquad
		\subfloat[$N=7,\ \nphy = 257^2$]{
			\includegraphics[width=0.27\linewidth]{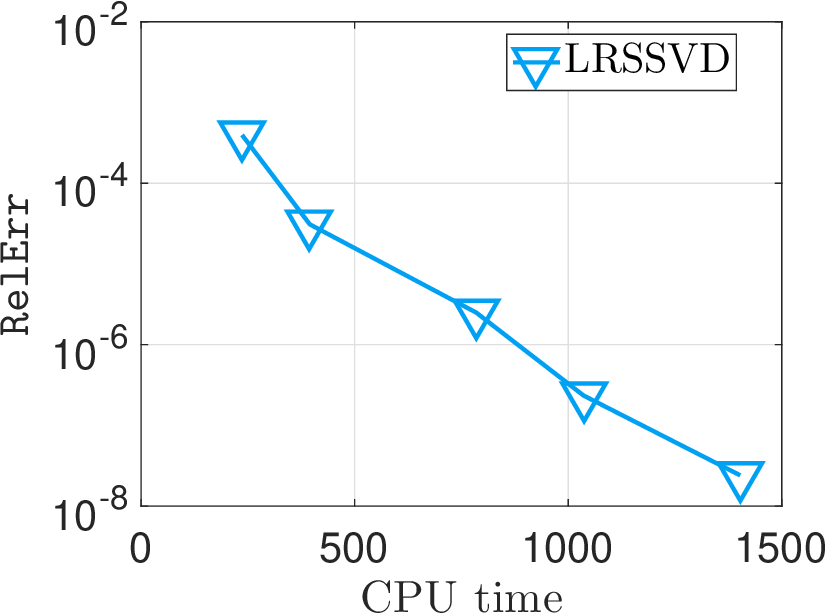}
		}\qquad
		\subfloat[$N=10,\ \nphy = 257^2$]{
			\includegraphics[width=0.27\linewidth]{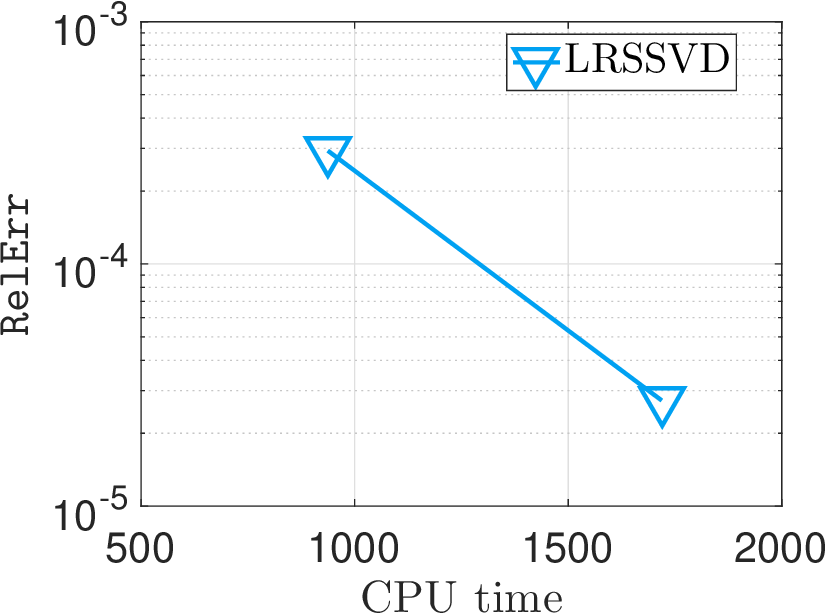}
		}
	\end{center}
	\caption{Relative errors with respect to CPU time in seconds for the LRSSVD approach.}\label{fig:test2e}
\end{figure}

Figure~\ref{fig:test2e} shows the CPU time in seconds for the LRSSVD approach. Note that for the case of $N=10$ and $\nphy=257^2$, the computation runs out of memory when the tolerance $\tol$ is set to $10^{-5}$, $10^{-6}$, and $10^{-7}$. From Figure~\ref{fig:test2c} and~\ref{fig:test2e}, it is clear that the AltLRP approach requires significantly less CPU time than the LRSSVD approach to achieve a given accuracy. In addition, if the required relative error is above $10^{-5}$, the SimLRP approach can also be applied.

\section{Conclusion}\label{sec:concusion}

In this work, we develop both a simultaneous low-rank projection (SimLRP) approach and an alternating low-rank projection (AltLRP) approach to compute the low-rank approximation of the solution for PDEs with random inputs. To identify the rank required for a desired accuracy, we propose a systematic strategy grounded in theoretical analysis, exploring the singular values of quasimatrices associated with truncated generalized polynomial chaos (gPC) expansions of the solution. Both the proposed methods exhibit enhanced computational efficiency compared to the conventional stochastic Galerkin method, particularly when dealing with numerous physical degrees of freedom. Numerical results  reveal that SimLRP excels in addressing low-dimensional problems, whereas AltLRP  demonstrates superior performance for moderate-dimensional problems. In the future, our focus will be on developing effective algorithms tailored for high-dimensional problems utilizing both AltLRP and SimLRP.

  \bibliographystyle{elsarticle-num} 
  \bibliography{ref}

\begin{thebibliography}{10}
\expandafter\ifx\csname url\endcsname\relax
  \def\url#1{\texttt{#1}}\fi
\expandafter\ifx\csname urlprefix\endcsname\relax\def\urlprefix{URL }\fi
\expandafter\ifx\csname href\endcsname\relax
  \def\href#1#2{#2} \def\path#1{#1}\fi

\bibitem{Xiu2002modeling}
D.~Xiu, G.~E. Karniadakis, Modeling uncertainty in steady state diffusion
  problems via generalized polynomial chaos, Computer Methods in Applied
  Mechanics and Engineering 191~(43) (2002) 4927--4948.
\newblock \href {https://doi.org/10.1016/S0045-7825(02)00421-8}
  {\path{doi:10.1016/S0045-7825(02)00421-8}}.

\bibitem{xiu2004two}
D.~Xiu, D.~M. Tartakovsky, A two-scale nonperturbative approach to uncertainty
  analysis of diffusion in random composites, Multiscale Modeling \& Simulation
  2~(4) (2004) 662--674.
\newblock \href {https://doi.org/10.1137/03060268X}
  {\path{doi:10.1137/03060268X}}.

\bibitem{Elman2005}
H.~C. Elman, O.~G. Ernst, D.~P. O’Leary, M.~Stewart, Efficient iterative
  algorithms for the stochastic finite element method with application to
  acoustic scattering, Computer Methods in Applied Mechanics and Engineering
  194 (2005) 1037--1055.
\newblock \href {https://doi.org/10.1016/j.cma.2004.06.028}
  {\path{doi:10.1016/j.cma.2004.06.028}}.

\bibitem{David2009}
O.~Ernst, C.~Powell, D.~Silvester, E.~Ullmann, Efficient solvers for a linear
  stochastic {Galerkin} mixed formulation of diffusion problems with random
  data, SIAM Journal on Scientific Computing 31~(2) (2009) 1424--1447.
\newblock \href {https://doi.org/10.1137/070705817}
  {\path{doi:10.1137/070705817}}.

\bibitem{Feng2015}
X.~Feng, J.~Lin, C.~Lorton, An efficient numerical method for acoustic wave
  scattering in random media, SIAM/ASA Journal on Uncertainty Quantification
  3~(1) (2015) 790--822.
\newblock \href {https://doi.org/10.1137/140958232}
  {\path{doi:10.1137/140958232}}.

\bibitem{Ghanem2003}
R.~G. Ghanem, P.~D. Spanos, Stochastic finite elements: a spectral approach,
  Courier Corporation, 2003.

\bibitem{Xiu2010}
D.~Xiu, Numerical methods for stochastic computations: a spectral method
  approach, Princeton University Press, 2010.

\bibitem{bedelman16}
B.~Soused\'ik, H.~Elman, {Stochastic Galerkin} methods for the steady-state
  {Navier-Stokes} equations, Journal of Computational Physics 316 (2016)
  435--452.
\newblock \href {https://doi.org/10.1016/j.jcp.2016.04.013}
  {\path{doi:10.1016/j.jcp.2016.04.013}}.

\bibitem{David2016}
A.~Bespalov, D.~Silvester, Efficient adaptive stochastic {Galerkin} methods for
  parametric operator equations, SIAM Journal on Scientific Computing 38~(4)
  (2016) A2118--A2140.
\newblock \href {https://doi.org/10.1137/15M1027048}
  {\path{doi:10.1137/15M1027048}}.

\bibitem{jin2016well}
S.~Jin, D.~Xiu, X.~Zhu, A well-balanced stochastic {Galerkin} method for scalar
  hyperbolic balance laws with random inputs, Journal of scientific computing
  67 (2016) 1198--1218.
\newblock \href {https://doi.org/10.1007/s10915-015-0124-2}
  {\path{doi:10.1007/s10915-015-0124-2}}.

\bibitem{Xiu2002wiener}
D.~Xiu, G.~E. Karniadakis, The {Wiener-Askey} polynomial chaos for stochastic
  differential equations, SIAM journal on scientific computing 24~(2) (2002)
  619--644.
\newblock \href {https://doi.org/10.1137/S1064827501387826}
  {\path{doi:10.1137/S1064827501387826}}.

\bibitem{Elman2014}
H.~C. Elman, D.~J. Silvester, A.~J. Wathen, Finite elements and fast iterative
  solvers: with applications in incompressible fluid dynamics, Oxford
  University Press (UK), 2014.

\bibitem{Pellissetti2000}
M.~F. Pellissetti, R.~G. Ghanem, Iterative solution of systems of linear
  equations arising in the context of stochastic finite elements, Advances in
  Engineering Software 31~(8) (2000) 607--616.
\newblock \href {https://doi.org/10.1016/S0965-9978(00)00034-X}
  {\path{doi:10.1016/S0965-9978(00)00034-X}}.

\bibitem{ullmann07}
M.~Eiermann, O.~G. Ernst, E.~Ullmann, Computational aspects of the stochastic
  finite element method, Computing and Visualization in Science 10~(1) (2007)
  3--15.
\newblock \href {https://doi.org/10.1007/s00791-006-0047-4}
  {\path{doi:10.1007/s00791-006-0047-4}}.

\bibitem{Powell2009}
C.~E. Powell, H.~C. Elman, Block-diagonal preconditioning for spectral
  stochastic finite element systems, IMA Journal of Numerical Analysis 29~(2)
  (2009) 350--375.
\newblock \href {https://doi.org/10.1093/imanum/drn014}
  {\path{doi:10.1093/imanum/drn014}}.

\bibitem{matthies2012solving}
H.~G. Matthies, E.~Zander, Solving stochastic systems with low-rank tensor
  compression, Linear Algebra and its Applications 436~(10) (2012) 3819--3838.
\newblock \href {https://doi.org/10.1016/j.laa.2011.04.017}
  {\path{doi:10.1016/j.laa.2011.04.017}}.

\bibitem{doostan2009least}
A.~Doostan, G.~Iaccarino, A least-squares approximation of partial differential
  equations with high-dimensional random inputs, Journal of computational
  physics 228~(12) (2009) 4332--4345.
\newblock \href {https://doi.org/10.1016/j.jcp.2009.03.006}
  {\path{doi:10.1016/j.jcp.2009.03.006}}.

\bibitem{powell2015}
C.~E. Powell, D.~Silvester, V.~Simoncini, An efficient reduced basis solver for
  stochastic {Galerkin} matrix equations, SIAM Journal on Scientific Computing
  39~(1) (2017) A141--A163.
\newblock \href {https://doi.org/10.1137/15M1032399}
  {\path{doi:10.1137/15M1032399}}.

\bibitem{LeeElman16}
K.~Lee, H.~C. Elman, A preconditioned low-rank projection method with a
  rank-reduction scheme for stochastic partial differential equations, SIAM
  Journal on Scientific Computing 39~(5) (2017) S828--S850.
\newblock \href {https://doi.org/10.1137/16M1075582}
  {\path{doi:10.1137/16M1075582}}.

\bibitem{elman2018low}
H.~C. Elman, T.~Su, A low-rank multigrid method for the stochastic steady-state
  diffusion problem, SIAM Journal on Matrix Analysis and Applications 39~(1)
  (2018) 492--509.
\newblock \href {https://doi.org/10.1137/17M1125170}
  {\path{doi:10.1137/17M1125170}}.

\bibitem{lee2019low}
K.~Lee, H.~C. Elman, B.~Sousedik, A low-rank solver for the {Navier--Stokes}
  equations with uncertain viscosity, SIAM/ASA Journal on Uncertainty
  Quantification 7~(4) (2019) 1275--1300.
\newblock \href {https://doi.org/10.1137/17M1151912}
  {\path{doi:10.1137/17M1151912}}.

\bibitem{wang2024reduced}
G.~Wang, Q.~Liao, Reduced basis stochastic {Galerkin} methods for partial
  differential equations with random inputs, Applied Mathematics and
  Computation 463 (2024) 128375.
\newblock \href {https://doi.org/10.1016/j.amc.2023.128375}
  {\path{doi:10.1016/j.amc.2023.128375}}.

\bibitem{benner2015low}
P.~Benner, A.~Onwunta, M.~Stoll, Low-rank solution of unsteady diffusion
  equations with stochastic coefficients, SIAM/ASA Journal on Uncertainty
  Quantification 3~(1) (2015) 622--649.
\newblock \href {https://doi.org/10.1137/130937251}
  {\path{doi:10.1137/130937251}}.

\bibitem{battles2005numerical}
Z.~Battles, Numerical linear algebra for continuous functions, Ph.D. thesis,
  University of Oxford (2005).

\bibitem{Townsend2014computing}
A.~Townsend, Computing with functions in two dimensions, Ph.D. thesis,
  University of Oxford (2014).

\bibitem{townsend2015continuous}
A.~Townsend, L.~N. Trefethen, Continuous analogues of matrix factorizations,
  Proceedings of the Royal Society A: Mathematical, Physical and Engineering
  Sciences 471~(2173) (2015) 20140585.
\newblock \href {https://doi.org/10.1098/rspa.2014.0585}
  {\path{doi:10.1098/rspa.2014.0585}}.

\bibitem{espig2014efficient}
M.~Espig, W.~Hackbusch, A.~Litvinenko, H.~G. Matthies, P.~W{\"a}hnert,
  Efficient low-rank approximation of the stochastic {Galerkin} matrix in
  tensor formats, Computers \& Mathematics with Applications 67~(4) (2014)
  818--829.
\newblock \href {https://doi.org/10.1016/j.camwa.2012.10.008}
  {\path{doi:10.1016/j.camwa.2012.10.008}}.

\bibitem{dolgov2015polynomial}
S.~Dolgov, B.~N. Khoromskij, A.~Litvinenko, H.~G. Matthies, Polynomial chaos
  expansion of random coefficients and the solution of stochastic partial
  differential equations in the tensor train format, SIAM/ASA Journal on
  Uncertainty Quantification 3~(1) (2015) 1109--1135.
\newblock \href {https://doi.org/10.1137/140972536}
  {\path{doi:10.1137/140972536}}.

\bibitem{Elman2007}
H.~Elman, D.~Furnival, Solving the stochastic steady-state diffusion problem
  using multigrid, IMA Journal of Numerical Analysis 27~(4) (2007) 675--688.
\newblock \href {https://doi.org/10.1093/imanum/drm006}
  {\path{doi:10.1093/imanum/drm006}}.

\bibitem{Babuska2004}
I.~Babu\v{s}ka, R.~l. Tempone, G.~E. Zouraris, Galerkin finite element
  approximations of stochastic elliptic partial differential equations, SIAM
  Journal on Numerical Analysis 42~(2) (2004) 800--825.
\newblock \href {https://doi.org/10.1137/S0036142902418680}
  {\path{doi:10.1137/S0036142902418680}}.

\bibitem{trefethen2022numerical}
L.~N. Trefethen, D.~Bau, Numerical linear algebra, Vol. 181, SIAM, 2022.

\bibitem{Caflisch1998}
R.~E. Caflisch, {Monte Carlo and quasi-Monte Carlo} methods, Acta numerica 7
  (1998) 1--49.
\newblock \href {https://doi.org/10.1017/S0962492900002804}
  {\path{doi:10.1017/S0962492900002804}}.

\bibitem{hoffman2003variation}
A.~J. Hoffman, H.~W. Wielandt, The variation of the spectrum of a normal
  matrix, in: Selected Papers Of Alan J Hoffman: With Commentary, World
  Scientific, 2003, pp. 118--120.

\bibitem{bhatia2007perturbation}
R.~Bhatia, Perturbation bounds for matrix eigenvalues, SIAM, 2007.

\bibitem{li2006eigenvalue}
W.~Li, J.-x. Chen, The eigenvalue perturbation bound for arbitrary matrices,
  Journal of Computational Mathematics (2006) 141--148.

\bibitem{Chen2019resfree}
Y.~Chen, J.~Jiang, A.~Narayan, A robust error estimator and a residual-free
  error indicator for reduced basis methods, Computers \& Mathematics with
  Applications 77~(7) (2019) 1963--1979.
\newblock \href {https://doi.org/10.1016/j.camwa.2018.11.032}
  {\path{doi:10.1016/j.camwa.2018.11.032}}.

\bibitem{jain2013low}
P.~Jain, P.~Netrapalli, S.~Sanghavi, Low-rank matrix completion using
  alternating minimization, in: Proceedings of the forty-fifth annual ACM
  symposium on Theory of computing, 2013, pp. 665--674.
\newblock \href {https://doi.org/10.1145/2488608.2488693}
  {\path{doi:10.1145/2488608.2488693}}.

\bibitem{lee2022enhanced}
K.~Lee, H.~C. Elman, C.~E. Powell, D.~Lee, Enhanced alternating energy
  minimization methods for stochastic {Galerkin} matrix equations, BIT
  Numerical Mathematics 62~(3) (2022) 965--994.
\newblock \href {https://doi.org/10.1007/s10543-021-00903-x}
  {\path{doi:10.1007/s10543-021-00903-x}}.

\bibitem{ifiss-siamreview}
H.~Elman, A.~Ramage, D.~Silvester, {IFISS}: A computational laboratory for
  investigating incompressible flow problems, SIAM Review 56 (2014) 261--273.
\newblock \href {https://doi.org/10.1137/120891393}
  {\path{doi:10.1137/120891393}}.

\bibitem{Berenger94}
J.~P. Berenger, A perfectly matched layer for the absorption of electromagnetic
  waves, Journal of Computational Physics 114~(2) (1994) 185--200.
\newblock \href {https://doi.org/10.1006/jcph.1994.1159}
  {\path{doi:10.1006/jcph.1994.1159}}.

\end{thebibliography}

\end{document}